\magnification=\magstep1
\input amstex
\documentstyle{amsppt}
\def \R {\Bbb R}
\def \C {\Bbb C}

\def \Z {\Bbb Z}

\def \T {\Bbb T}
\def \F {\Bbb F}
\def \H {\Bbb H}
\def \O {\Bbb O}
\def \A {\Cal A}

\def \E {\Cal E}

\def \G {\Cal G}

\def \W {\Cal W}
\def \g {\frak g}
\def \h {\frak h}
\def \k {\frak k}

\def \a {\frak a}

\def\gl {\frak g\frak l}

\def\p {\frak p}
\def \n {\frak n}

\def \u {\frak u}
\def \v {\frak v}
\def \w {\frak w}
\def \z {\frak z}

\def \RE {\Re\text{\rm e}\,}
\def \IM {\Im\text{\rm m}\,}
\def\one{{\text{\bf 1}}}
\def \al {\alpha}
\def \la {\lambda}
\def \ph {\varphi}
\def \del {\delta}
\def \eps {\varepsilon}
\def\om {\omega}
\def \lan {\langle}
\def \ran {\rangle}
\def \proof {\demo {Proof}}
\def \endproof {\qed\enddemo}

\def \trans{\,{}^t\!}
\def \half{\frac12}
\def \inv{^{-1}}

\def \Dot {\overset . \to}

\def \vol {\text{\rm vol\,}}
\def \dim {\text{\rm dim\,}}
\def \span {\text{\rm span\,}}
\def \rad {\text{\rm rad\,}}

\def \Cliff{\text{\rm Cliff}}
\def\id {\text{\rm id}}

\def\Pin {\text{\rm Pin}}
\def\Ad {\text{\rm Ad}}

\topmatter
\title A unified approach to compact symmetric spaces of rank one
\endtitle
\rightheadtext{Rank-one symmetric spaces}
\leftheadtext{A. Kor\'anyi and F. Ricci}
\author Adam Kor\'anyi
 and Fulvio Ricci
\endauthor
\address Department of Mathematics and Computer Science, Lehman
College, City University of New York, Bronx, NY 10468-1589, USA
\endaddress
\address Scuola Normale Superiore, Piazza dei Cavalieri
7, 56126 Pisa, Italy \endaddress
\thanks Research partially supported by the National Science Foundation, a
PSC-CUNY grant, Centro de Giorgi, Pisa and GNAMPA-INdAM, Rome.
\endthanks
\abstract A relatively simple algebraic framework is given, in which all the compact symmetric spaces can be described and handled without distinguishing cases. We also give some applications and further results.
\endabstract
\endtopmatter

\document

\head Introduction\endhead
\vskip.3 cm

This article deals with compact Riemannian globally symmetric spaces of rank one. These spaces are very important for geometry, and of course they are well known. They are the spheres and the projective spaces over $\R,\C,\H,\O$ (in the last case only of dimension $\le2$).

There are two standard ways of describing them and make computations on them. One is to use classification; this is the approach, e.g., in \cite{Be}. The projective spaces over $\R,\C,\H$ are easily handled together, but the case of $\O$ is different because of the non-associativity of the octonions. For the $\O$-case one has to refer, as does \cite{Be}, to the rather complicated articles of Freudenthal (cf., e.g., \cite{Ba, S}).

The other way is is to use the general theory of symmetric spaces \cite{Hel1}. Here there is no need to distinguish cases, but one has to use the large machinery of semisimple Lie group theory, which was designed for much more general situations and is rather unwieldy when applied to the special case of spaces of rank one.

In this article we set up a new framework which does not make explicit use of the octonions and applies to spaces of rank one in an essentially unified way and makes it easy to work on them, including the $\O$-case. This is done in the first six sections.

The basic notion is that of a {\it $C$-module with $J^2$-condition} (briefly: $J^2C$-module). A $C$-module is the same thing as the classical notion of composition of quadratic forms \cite{Hus}, and is closely related to the notion of a Clifford module. The $J^2$-condition specifies a subclass; it already played a fundamental r\^ole in \cite{CDKR2}. As we have discovered since, its basic idea goes back to \cite{Hei}.

For non-compact symmetric spaces, a program similar to the present one was carried out in \cite {CDKR2} and our article can be regarded as a continuation of \cite{CDKR2}, although our starting point is somewhat different. In \cite{CDKR2} the basic objects are $H$-type Lie algebras; they form a category essentially equivalent to that of Clifford modules.

Some of the arguments in the introductory sections of this article are classical or are reformulations of proofs in \cite{CDKR2}, but we prove everything we need. Some other facts, which we do not need but which are of interest for understanding the full picture, are only described without proof  in Section 3. This section also contains the classification of $J^2C$-modules, with proofs. This, together with the known classification of all symmetric spaces \cite{Hel1, W}, could be used to quickly show that our theory covers all the compact rank-one symmetric spaces. However, we will give a classification-free proof of this fact in Section 8.

Sections 1 and 2 contain the basic properties of $J^2C$-modules. In Section 4 we introduce the space $W=C\oplus V$, with $V$ a $J^2C$-module. $W$ is a weaker substitute for a $J^2C$-module, but  a good notion of $C$-line can still be defined in it.

With these tools, we construct $CPW$ as a compactification of $W$ by adjoining a point at infinity for every family of parallel $C$-lines in $W$. 
In Section 5 we describe the topology and the differentiable structure on $CPW$, and in Section 6 its metric and its isometry group $U$, proving that $CPW$ actually is a compact rank-one symmetric space.
This ends the construction of the symmetric spaces. 

In the subsequent part of the paper we prove various facts about the geometry of these spaces and the relevant transformation groups on them.

In Section 7 we illustrate the usefulness of our construction by reproving in simple ways some known properties of these spaces. As mentioned already, we prove in Section 8 that they are all the compact rank-one symmetric spaces, without appealing to classification.

In Section 9 we analyze the structure of the group $GL(W,C)$ of invertible linear transformations of $W$ which preserve $C$-lines. This is preliminary to Section 10, where we introduce the group $\G$ of collineations of $CPW$. $\G$ is a transformation group properly containing the isometry group $U$ and $GL(W,C)$, and whose elements are characterized by the property of preserving projective $C$-lines (i.e. closures in $CPW$ of the $C$-lines of $W$, or images of these under isometries). We prove that collineations form a semisimple Lie group and that they act conformally on each projective $C$-line.

In Section 11 we prove that the collineation groups $\G$ are characterized as the semisimple parts of automorphism groups of irreducible symmetric cones. We do so by constructing a representation of $\G$ and identifying the cone in a real form of the representation space. The induced projective action of $\G$ on the space of extremal lines of the cone provides an identification of this space with $CPW$.

In the Appendix at the end of the paper we show how the general non-compact symmetric space can be realized as the unit ball in $W$ with a different metric. This amounts to redoing \cite{CDKR2} from a different starting point.

We wish to point out that our construction results in an actual unification in the geometric as well as the algebraic sense. The spheres are included as the extreme case $V=0$ and $C$ arbitrary. In such a case the whole $CPW$ consists of a single projective $C$-line. The other extreme case, $C=\R$ and $V$ arbitrary, gives real projective space. In both extreme cases the unit ball of $W$ is real hyperbolic space; in the first case it appears as the Poincar\'e model, in the second as the Klein model. In the cases in between, it agrees with the models used in \cite{M}.

\vskip.5cm

\head Table of contents\endhead

\roster
\item"1." $C$-modules with the $J^2$-condition
\item"2." Automorphisms of $J^2C$-modules
\item"3." Classification and some background information
\item"4." The space $W=C\oplus V$
\item"5." The construction of $CPW$
\item"6." $CPW$ as a compact symmetric space
\item"7." Some applications
\item"8." Every compact rank-one space is a $CPW$
\item"9." The group $GL(W,C)$
\item"10." The group of collineations
\item"11." Compact rank-one spaces and symmetric cones
\item"A." Appendix. The non-compact symmetric spaces
\endroster

\vskip.5cm

\head 1. $C$-modules with the $J^2$-condition \endhead

\vskip.3cm

Let $C$, $V$ be finite-dimensional real Euclidean spaces over $\R$. We write the inner products as $\lan\cdot ,\cdot \ran$ and the norms as $|\cdot|$. We assume that $C$ has a distinguished unit element, denoted by $\one$. 

A {\it $C$-module structure} on $V$ is a bilinear map $J:C\times V\rightarrow V$, such that
$$
\align
J(\one,v)&=v\qquad (\forall\,v\in V)\ ,\tag1.1\\
\big|J(\zeta,v)\big|&=|\zeta||v|\qquad (\forall\,\zeta\in C\,,\, v\in V)\ .\tag1.2
\endalign
$$

Instead of $J(\zeta,v)$ we will also use the notation $J_\zeta v$ and, more frequently, $\zeta v$. We also write $Cv=\{\zeta v:\zeta\in C\}$ for $v\in V$.

We say that the $C$-module $V$ {\it satisfies the $J^2$-condition} (or, briefly, is a $J^2C$-module) if 
$$
C(Cv)=Cv\qquad (\forall\,v\in V)\ .\tag1.3
$$

We do not exclude the ``trivial'' cases where $C=\R\one$ or $V=0$. Notice that in such cases the $J^2$-condition is trivially satisfied and most of the next definitions are vacuous.

We denote by $C'$ the orthogonal complement of $\R\one$ in $C$. If $\zeta=a\one+z$, with $a\in\R$ and $z\in C'$, we set $\bar\zeta=a\one-z$, $a=\RE\zeta$, $z=\IM\zeta$.

Polarizing (1.2) in both $\zeta$ and $v$, we have
$$
\lan \zeta u,\eta v\ran+\lan \eta u,\zeta v\ran=2\lan \zeta,\eta \ran\lan  u, v\ran\ ,\tag1.4
$$

Taking $\eta=\one$, $\zeta=z\in C'$, we see that $J_z$ is skew-symmetric. For general $\zeta\in C$, this implies that
$$
J_{\bar\zeta}J_\zeta=|\zeta|^2\id\ .\tag1.5
$$

For $z\in C'$, (1.5) gives $J_z^2=-|z|^2\id$. From this it is clear that the action of $C'$ on $V$ via $J$ extends to an action of the Clifford algebra $\Cliff(C')$. (We recall that $\Cliff(C')$ is the associative algebra generated by $C'$ and a unit element $e$ subject to the relations $z^2=-|z|^2e$. In other words, if $\{z_1,\dots,z_m\}$ is an orthonormal basis of $C'$, then $\Cliff(C')$ is generated by these elements and $e$, subject to the relations $z_iz_j+z_jz_i=-2\del_{ij}e$.)

So $V$ is automatically a {\it Clifford module}, meaning a representation of the associative algebra $\Cliff(C')$. (More exactly we should say an ``orthogonal'' Clifford module, to take into account the added condition $|zv|=|z||v|$ for all $z\in C'$, $v\in V$.)

For any $\al\in\Cliff(C')$, we still write its action on $v\in V$ as $\al v$.

By (1.3), for every $v\ne0$ in a $J^2C$-module $V$, $Cv$ is a $\Cliff(C')$ submodule, necessarily irreducible. From (1.4), it follows that the orthogonal complement of a submodule is again a submodule. So we can inductively construct {\it orthonormal $C$-bases}, i.e. orthonormal sets $\{v_1,\dots,v_n\}$ such that $V$ is the orthogonal sum of $Cv_1, \dots, Cv_n$. If $V=Cv$ for some $v\ne0$, we say that the $J^2C$-module is irreducible.

Another important property of $J^2C$-modules is that any non-zero element $v\in V$ determines a multiplication law on $C$, denoted by $\cdot_v$ and given by
$$
(\zeta\cdot_v\eta)v=\zeta\eta v\ .
$$

Obviously, $\cdot_{tv}=\cdot_v$ for $t\in\R^*$.

\vskip.3cm

For the remainder of this section, we assume that $V$ is a non-trivial $J^2C$-module (but we allow $C=\R\one$).

\proclaim{Proposition 1.1} Under the multiplication $\cdot_v$, $C$ is a normed division algebra.
\endproclaim

\proof The product is well defined because the map $\zeta\mapsto\zeta v$ from $C$ to $Cv$ is a bijection. It is obvious that the product is bilinear and has $\one$ as its identity element. Also, $|\zeta\cdot_v\eta|=|\zeta||\eta|$ by (1.2). To show that $C$ is a division algebra, we must solve the equations
$$
\xi\cdot_v\eta=\zeta\ ,\qquad \eta\cdot_v\xi'=\zeta
$$ 
in $\xi,\xi'$, when $\zeta$ and $\eta\ne0$ are given. For the first equation, we have to solve
$$
\xi\eta v=\zeta v\ .
$$

This can be done, because the map $\la\mapsto \la(\eta v)$ is a bijection of $C$ onto $C(\eta v)\subset C(Cv)=Cv$. Since both $Cv$ and $C(\eta v)$ have the same dimension as $C$, we have $C(\eta v)=Cv$. 

To solve the equation
$$
\eta \xi' v=\zeta v\ ,
$$
we multiply both sides by $\bar\eta$ on the left and use (1.5) to obtain that
$$
\xi' v=|\eta|^{-2}\bar\eta\zeta v\ ,
$$
i.e. $\xi'=|\eta|^{-2}\bar\eta\cdot_v\zeta$.
\endproof

\vskip.2cm

\subhead Remark \endsubhead
As we shall see next in detail, $\cdot_v$ can actually depend on $v$, and this dependence is equivalent to lack of associativity. Nevertheless, some expressions are independent of $v$ in general. We list a few.
\roster
\item"1." The proof of Proposition 1.1 shows that $\zeta\inv=|\zeta|^{-2}\bar\zeta$ when $\zeta\ne0$, with respect to any multiplication $\cdot_v$. It follows that the value of any other rational expression in $\zeta$ with real coefficients is independent of $v$ and coincides with its value with respect to Clifford multiplication.
\item"2." For any $\zeta,\eta\in C$, the product $\zeta\cdot_v(\eta\cdot_v\zeta)$ equals the Clifford algebra product $\zeta\eta\zeta$ and therefore it does not depend on $v$. To see this, it is convenient to consider separately the case $\eta\in\R\one+\R\IM\zeta$ and the case $\eta\perp(\R\one+\R\IM\zeta)$. The first case is trivial. In the second case we set $\zeta=a\one+z$. Since $\eta\in C'$ and $\eta\perp z$,
$$
\align
\zeta\eta\zeta &=(a\one+z)\eta(a\one+z)\\
&=a^2\eta +a(z\eta+\eta z)+z\eta z \\
&=a^2\eta +|z|^2\eta \ ,
\endalign
$$
which is in $C$.
\item"3."  For any $v$, $\RE(\zeta\cdot_v\eta)=\lan \zeta,\bar\eta\ran$. In fact, both expressions are symmetric, bilinear, and they agree when $\zeta=\eta$.
\endroster

\vskip.3cm

\proclaim{Proposition 1.2} Any division subalgebra of $(C,\cdot_v)$ generated by two elements is associative.
\endproclaim

This is true in any normed division algebra. An easy classification-free proof is in \cite{FK, pp. 82-83}. It is also well known and easy (cf. e.g. \cite{FK}) that the only normed division algebras are $\R,\C,\H,\O$. But, for self-sufficiency of this article, here is a proof of the Proposition, making use of the present setup.

\proof We look at the intersection of the subalgebra with $C'$. If it has dimension 0 or 1, the conclusion is trivial. If it has dimension at least 2, we can assume that the two generators $z_1,z_2$ are both in $C'$ and orthonormal. Let $z_3=z_1\cdot_vz_2$. Then $z_3\in C'$ and orthogonal to both $z_1$ and $z_2$. In fact,
$$
\align
\lan z_3,\one\ran|v|^2&=\lan z_1z_2v,v\ran=-\lan z_2v,z_1v\ran=-\lan z_2,z_1\ran|v|^2=0\ ,\\
\lan z_3,z_1\ran|v|^2&=\lan z_1z_2v,z_1v\ran=\lan z_2v,v\ran=\lan z_2,\one\ran|v|^2=0\ ,
\endalign
$$
and similarly $z_3\perp z_2$. Furthermore,
$$
(z_1\cdot_v z_3)v=z_1(z_1z_2v)=z_1^2z_2v=-z_2v\ ,
$$
so $z_1\cdot_v z_3=-z_2$. Similarly, $z_2\cdot_v z_3=z_1$. The Clifford relations imply that $z_i\cdot_v z_j=-z_j\cdot _v z_i$. This shows that $(C,\cdot_v)\cong \H$, which is associative.
\endproof

\proclaim{Corollary 1.3} If $\zeta,\eta,\lambda\in C$ belong to a division subalgebra of $(C,\cdot_v)$ generated by two elements, and $\la\ne0$, then 
$$
\zeta\cdot_v\eta=\zeta\cdot_{\la v}\eta\ .
$$
\endproclaim

\proof We have
$$
\align
(\zeta\cdot_{\la v}\eta)(\la v)&=\zeta\eta\la v=\zeta(\eta\cdot_v\la)v\\
&=\big(\zeta\cdot_v(\eta\cdot_v\la)\big)v=\big((\zeta\cdot_v\eta)\cdot_v\la\big)v\\
&=(\zeta\cdot_v\eta)(\la v)\ ,
\endalign
$$
where we have used associativity of $\cdot_v$ in the second line.
Then use the fact that the map $\tau\mapsto \tau(\la v)$ from $C$ to $C(\la v)$ is one-to-one.
\endproof

\proclaim{Proposition 1.4} If $(C,V)$ is not irreducible, then the multiplication $\cdot_v$ is independent of $v$.
\endproclaim

The proof we give is also contained in the proof of Theorem 1.1 of \cite{CDKR1}.

\proof Given $u,v$ non-zero elements of $V$, we prove that $\cdot_u=\cdot_v$. Assume first that $Cu+Cv$ is a direct sum. For $\zeta,\eta\in C$,
$$
(\zeta\cdot_{u+v}\eta)(u+v)=(\zeta\cdot_{u+v}\eta)u+(\zeta\cdot_{u+v}\eta)v\ ,
$$
but also
$$
(\zeta\cdot_{u+v}\eta)(u+v)=\zeta\eta(u+v)=\zeta\eta u+\zeta\eta v=(\zeta\cdot_u\eta)u+(\zeta\cdot_v\eta)v\ .
$$

Therefore $(\zeta\cdot_{u+v}\eta)u=(\zeta\cdot_u\eta)u$ and $(\zeta\cdot_{u+v}\eta)v=(\zeta\cdot_v\eta)v$, i.e. $\zeta\cdot_u\eta=\zeta\cdot_{u+v}\eta=\zeta\cdot_v\eta$.

Assume now that $Cu+Cv$ is not a direct sum. If $\zeta_0u=\eta_0v$ with $\zeta_0,\eta_0\ne0$, then $
v=\eta_0\inv\zeta_0u\in C(Cu)=Cu$, and so $Cv=Cu$. Since $(C,V)$ is not irreducible, there is $v'\ne0$ such that $v'\perp Cv$. Using the first part of the proof, $\cdot_v=\cdot_{v'}=\cdot_u$.
\endproof

\proclaim{Corollary 1.5} The product $\cdot_v$ is independent of $v$ if and only if $(C,\cdot_v)$ is associative for one (and hence for all) $v$. This is the case when $\dim C\le4$.
\endproclaim

\proof We have the general identity
$$
\zeta\cdot_v(\eta\cdot_v\la)=(\zeta\cdot_{\la v}\eta)\cdot_v\la\qquad(\forall\,\zeta,\eta,\la\in C)\ .\tag1.6
$$

Then, if $\cdot_v$ is independent of $v$, then $(C,\cdot_v)$ is clearly associative. Assume now that, for some $v\ne0$, $(C,\cdot_v)$ is associative. By Proposition 1.4, we can restrict ourselves to the irreducible case $V=Cv$.
So any other non-zero $u\in V$ is equal to $\la v$ for some non-zero $\la\in C$. By (1.6),
$$
(\zeta\cdot_v\eta)\cdot_v\la=\zeta\cdot_v(\eta\cdot_v\la)=(\zeta\cdot_u\eta)\cdot_v\la\ .
$$

Dividing by $\la$, we have $\zeta\cdot_v\eta=\zeta\cdot_u\eta$.

If $\dim C=1$, then $C=\R\one$ and there is nothing to say. If $\dim C=2$, taking $i\in C'$ with $i^2=-1$, it is immediate to see that $C\cong\C$. If $\dim C\ge3$, take orthonormal vectors $i,j\in C'$. Fix $v\in V$ a non-zero vector, and set $k=i\cdot_vj$. As shown in the proof of Proposition 1.2, $k$ is linearly independent of $\one,i,j$, so that $\dim C=4$ and $(C,\cdot_v)$ is generated by $i$ and $j$. By Proposition 1.2, $(C,\cdot_v)\cong\H$ is associative.
\endproof

We shall say in short that $C$ is associative, or non-associative, to distinguish between these two cases. 

\vskip.5cm

\head 2. Automorphisms of $J^2C$-modules \endhead

\vskip.3cm

Let $(C,V)$ be a $C$-module. An automorphism of $(C,V)$ is a pair $m=(\ph,\psi)$ of orthogonal maps $\ph:C\rightarrow C$, $\psi:V\rightarrow V$ such that the diagram
$$
\matrix C\times V&\overset J\to\longrightarrow& V\\
\ph\downarrow\ \ \ \ \downarrow \psi&&\ \ \ \downarrow\psi\\
C\times V&\overset J\to\longrightarrow&  V
\endmatrix
$$
is commutative, i.e.
$$
\psi(\zeta v)=\ph(\zeta)\psi(v)\ .\tag2.1
$$

 We write $M$  for the group of automorphisms of the $C$-module. The automorphism group of a general $C$-module is described in detail in \cite{R}. We shall look more closely at the specific properties of $M$ when $V$ satisfies the $J^2$-condition.

We denote by $M_1$ the subgroup of $M$ defined by the condition $\ph=\id$. $M_1$ is the automorphism group of the Clifford module associated to $(C,V)$.

\proclaim{Proposition 2.1} Let $(C,V)$ be a $J^2C$-module, and assume that $C$ is associative. Then $M_1$ is in one-to-one correspondence with the ordered orthonormal $C$-bases of $V$.
\endproclaim

\proof Let $\{u_1,\dots,u_n\}$ and $\{v_1,\dots,v_n\}$ be two o.n. $C$-bases. Define $\psi:V\rightarrow V$ by
$$
\psi\Big(\sum \zeta_ju_j\Big)=\sum \zeta_jv_j\ .
$$

This is clearly a well-defined orthogonal map. To prove that $(\id,\psi)\in M$, we must verify that 
$$
\psi(\eta v)=\eta\psi(v)\qquad (\forall\,\eta\in C,v\in V)\ .\tag2.2
$$ 

Since $\cdot_v$ is independent of $v$,
$$
\align
\psi\Big(\eta\sum \zeta_ju_j\Big)&=\psi\Big(\sum (\eta\cdot_{u_j}\zeta_j)u_j\Big)\\
&=\sum (\eta\cdot_{u_j}\zeta_j)v_j\\
&=\sum (\eta\cdot_{v_j}\zeta_j)v_j\\
&=\eta\sum \zeta_jv_j\ .
\endalign
$$

It is then obvious that $\psi$ is the only linear map satisfying (2.2) and mapping each $u_j$ into $v_j$.
\endproof

To construct automorphisms that are not trivial on $C$ we need a preliminary remark about conjugation in $\Cliff(C')$.

If $z$ is a unit element of $C'$ and $\eta\in C$, then the Clifford product $z\eta z\inv=-z\eta z$ is also in $C$ (see Remark 2 in Section 1). In fact, it equals $\eta$ if $\eta \in\R\one+\R z$, and $-\eta$ if $\eta\perp(\R\one+\R z)$ (in other words, $z\eta z\inv$ is the reflection of $\eta$ in the plane $\R\one+\R z$).

Let $\Pin(C')\subset\Cliff(C')$ be the multiplicative group generated by unit elements of $C'$. It follows that, if $\al\in\Pin(C')$, the map $\eta\mapsto \al \eta\al\inv$ is orthogonal on $C$ and is the identity on $\R\one$.

On a general $C$-module, for any  $\al\in\Pin(C')$,
$$
m_\al(\zeta,v)=(\al\zeta\al\inv,\al v)\ ,\tag2.3
$$
is in $M$. The $m_\al$ form a subgroup that we denote by $M_2$. $M_2$ acts on $C'$ as $SO(C')$ if $\dim C'$ is odd, and as $O(C')$ if $\dim C'$ is even (this follows because this action is generated by the reflections $\eta\mapsto z\eta z\inv$). In particular, $M_2$ acts transitively on spheres in $C'$. It follows that $M_2$ is isomorphic with $\Pin(C')$ modulo a finite subgroup. It is also true, even if of no great importance for us here, that $M_1\cap M_2$ is finite and that the group $M_1M_2$ has index at most 2 in $M$ (cf. \cite{R}).

\proclaim{Proposition 2.2} Let $(C,V)$ be a $J^2C$-module, and let $\Xi_V$ be the manifold of ordered orthonormal $C$-bases of $V$. Let also $S_{C'}$ be the unit sphere in $C'$. Then $M$ acts transitively on $S_{C'}\times \Xi_V$.
\endproclaim

\proof For any $m=(\ph,\psi)\in M$, $\ph(\one)=\one$, hence $\ph$ acts orthogonally on $C'$. It is also simple to verify that $\psi$ transforms orthonormal $C$-bases of $V$ into orthonormal $C$-bases. Therefore $M$ acts on $S_{C'}\times\Xi_V$. 

Given two elements $(z_0,\xi_0), (z'_0,\xi'_0)\in S_{C'}\times\Xi_V$, we want to find $m\in M$ such that $m(z'_0,\xi'_0)=(z_0,\xi_0)$. There is $\al\in\Pin(C')$ such that $\al z'_0\al\inv=z_0$. Then $\al$ transforms $\xi'_0$ into another o.n. basis $\xi'_1$. Then we need to find $m\in M$ such that $m(z_0,\xi'_1)=(z_0,\xi_0)$.

If $C$ is associative, Proposition 2.1 says that there is an element of $M_1$ that does the job. If $C$ is non-associative, we necessarily have $\dim C=\dim V\ge5$, by Proposition 1.4 and Corollary 1.5. In particular, $\Xi_V= S_V$, the unit sphere in $V$, since $(C,V)$ is irreducible. We then have two elements $v_1,v_2\in S_V$, and we want to find $\beta\in\Pin(C')$ such that $\beta z_0\beta\inv=z_0$ and $\beta v_1=v_2$.

Choose $v\in S_V$ orthogonal to $v_j$ and $z_0v_j$ for $j=1,2$. Since $V=Cv$, there are $z_1,z_2\in C$ such that $v_j=z_jv$, $j=1,2$. By (1.4), $v\perp v_j$ implies that $z_j\in C'$, and $v\perp z_0v_j$ implies that $z_j\perp z_0$. Let $\beta=-z_2z_1=z_2z_1\inv$. Then
$$
\align
\beta z_0\beta\inv&=z_2z_1z_0z_1z_2=-z_2z_0z_2=z_0\ ,\\
\beta v_1&=z_2z_1\inv z_1v=z_2v=v_2\ .\qquad\quad\qed
\endalign
$$
\enddemo

\proclaim{Corollary 2.3} $M$ acts transitively on the product $S_{C'}\times S_V$ of the two unit spheres in $C'$ and $V$, and on the product of unit spheres $S_{C'}\times S_{Cv_1}\times\cdots\times S_{Cv_n}$, if $\{v_1,\dots,v_n\}$ is an orthonormal $C$-basis of $V$.
\endproclaim

 Transitivity of $M$ on $S_{C'}\times S_V$ is {\it Kostant's double transitivity} for the $\Ad(M)$-action on the sum of root spaces $\g^\la+\g^{2\la}$ in semisimple Lie algebras.

\vskip.5cm

\head 3. Classification and some background information \endhead

\vskip.3cm

Theorem 3.1 below gives the classification of all $J^2C$-modules in terms of normed division algebras. This theorem will not be used at all in the rest of the article. 

\proclaim{Theorem 3.1} Every $J^2C$-module is isomorphic with one of the following:
\roster
\item"(i)" $C$ is any Euclidean space with a distinguished unit vector $\one$, $V=0$;
\item"(ii)" $C=\F$, a normed division algebra, $V=\F^n$ (with $n$ a positive integer if $\F$ is associative, and $n=1$ if $\F$ is non-associative, the norm on $\F^n$ being the usual $\ell^2$-norm) and $J$ is multiplication by elements of $\F$ from the left.
\endroster
\endproclaim

\proof These are clearly $J^2C$-modules. The $J^2$-condition is vacuous in (i) and trivial in (ii) when $\F$ is associative. If $\F$ is non-associative, then $Cv=V$ for any $v\ne0$, and the $J^2$-condition follows trivially.

To prove that our list is complete, assume that $(C,V)$ is a $J^2C$-module. If $V=0$, there is nothing to prove. If $V\ne0$, we must provide an isomorphism $(f,g)$ from some $(\F,\F^n)$ as above to $(C,V)$, i.e. a pair of orthogonal linear maps $f:\F\rightarrow C$, $g:\F^n\rightarrow V$ such that
$$
g\big(q(q_1,\dots,q_n)\big)=f(q)g(q_1,\dots,q_n)\qquad (\forall, q,q_1,\dots,q_n\in\F)\ .
$$

Fix an orthonormal $C$-basis $\{v_1,\dots,v_n\}$ of $V$ and set $\F=(C,\cdot_{v_1})$. Let $f:\F\rightarrow C$ be the identity map and define $g:\F^n\rightarrow V$ by
$$
g(\zeta_1,\dots,\zeta_n)=\sum\zeta_jv_j\ .
$$

If $\F$ is associative, $\cdot_{v_1}=\cdot_{v_j}$ for every $j$. Then (as in the proof of Proposition 2.1)
$$
\align
g\big(\eta\cdot_{v_1}(\zeta_1,\dots,\zeta_n)\big)&=\sum (\eta\cdot_{v_1}\zeta_j)v_j\\
&=\sum\eta\zeta_jv_j\\
&=\eta\sum\zeta_jv_j\\
&=f(\eta)g(\zeta_1,\dots,\zeta_n)\ .
\endalign
$$

This computation also works when $C$ is non-associative, due to the fact that $n=1$.
\endproof

\vskip.3cm

We conclude this section by mentioning the relations among $C$-modules, Clifford modules and $H$-type Lie algebras.

As pointed out in Section 1, every $C$-module extends in a natural way to a $\Cliff(C')$-module.
This construction can be reversed: let $C'$ be any finite dimensional (possibly trivial)
vector space with a scalar product and
$V$ a module over $\Cliff(C')$, endowed with a $\Pin(C')$-invariant scalar
product. Setting $C=\R \one\oplus C'$, there is a unique $C$-module structure
$(C,V)$ inducing the given Clifford module structure on $V$
\cite{K}. However, non-equivalent Clifford modules may induce isomorphic
$C$-modules. This is because $(C,V)$ and $(C,V')$ may be isomorphic without any isomorphism being the identity on $C$ (cf. \cite{KR}).

The main facts about Clifford modules can be derived from
\cite{Hus}. If the dimension $d$ of $C$ is not divisible by 4, $\Cliff(C')$ has, up to
equivalence, only one irreducible module $V_0$, and every
other module can be realized as $V_k\sim V_0\otimes_\F \F^k$, where
$\F=\R,\C,\H$, depending on the congruence class of $d$ mod.\,8. 
If $d$ is divisible by 4, then $\Cliff(C')$ has two inequivalent irreducible
modules, $V_1,V_2$, and any module $V$ can be realized as
$V_{kh}=(V_1\otimes_\F\F^k)\oplus(V_2\otimes_\F\F^h)$, with $\F=\R$
if $d$ is divisible by 8, and $\F=\H$ otherwise. If $k\ne h$, $V_{kh}$
and $V_{hk}$ are non-equivalent as Clifford modules, but the induced
$C$-modules are isomorphic.

The $J^2$-condition makes sense on Clifford modules, and it can be restated as $\Cliff(C')v=Cv$ for every $v\in V$. This can only occur when $\Cliff(C')$ has a module (necessarily irreducible) of the same dimension as $C$, i.e. when $d=1,2,4,8$. For these values of $d$, the $J^2$-condition holds for every module if $d=1,2$, only for the ``isotypic'' $C$-modules $V_{k0}=V_{0k}$ if $d=4$, and only for the irreducible $C$-module if $d=8$ \cite{CDKR1}.

\vskip.3cm

There is a third category which must be mentioned in this context: the $H$-type Lie algebras.
An $H$-type Lie algebra (the notion of which is the basic object in \cite{CDKR1,2}) is a Lie algebra $\n$ with a positive definite inner product. It is assumed that $\z$ is the center of $\n$, $\n=\v+\z$ is an orthogonal direct sum and $[\v,\v]\subset\z$. Furthermore, for all $z\in\z$, the map $J_z:\v\rightarrow\v$, defined by
$$
\lan J_zv,u\ran=\big\lan[v,u],z\big\ran\ ,\tag3.1
$$
has the property
$$
J_z^2=-|z|^2\id_\v\ .\tag3.2
$$

It is immediate that, given a $C$-module $(C,V)$, we obtain an $H$-type Lie algebra by taking $\n$ as the Euclidean direct sum of $\v=V$ and $\z=C'$ and defining the Lie bracket through (3.1). The converse is also easy to verify, as well as that isomorphic $C$-modules produce isomorphic Lie algebras, and {\it vice versa} \cite{K, KR}.

\vskip.5cm

\head 4. The space $W=C\oplus V$ \endhead

\vskip.3cm

Given a $J^2C$-module $(C,V)$, let $W=C\oplus V$ be the direct sum of Euclidean spaces. We introduce the following equivalence relation on $W\setminus\{0\}$:
\vskip.2cm
\roster
\item"(i)" $(0,u)\sim(0,v)$ if $u\in Cv$;
\item"(ii)" if $\zeta\ne0$, $(\eta,u)\sim(\zeta, v)$ if $\eta\ne0$ and $\eta\inv u=\zeta\inv v$.
\endroster

\vskip.2cm

Notice that the $J^2$-condition is required in (i) to prove transitivity. Given an element $w=(\zeta,v)\ne0$ in $W$, we denote its equivalence class together with point 0 by $Cw$ or $C(\zeta,v)$. We must pay attention, however, to the fact that the notation $\eta w$ does not make any sense for an individual $\eta\in C$. This is related to the fact that in general $W$ cannot be given a $C$-module structure in a natural way (unless $C$ is associative). We only have that, for $v\ne0$, $C(\zeta,v)=\big\{(\la\cdot_v\zeta,\la v):\la\in C\big\} =\big\{(\la,\la\zeta\inv v):\la\in C\big\}$.

We call {\it $C$-line}, or {\it affine $C$-line}, any translate $w_0+Cw$ of $Cw$, with $w\ne0$. The $C$-lines through 0 form a closed subset of the (Grassmanian) manifold of real $d$-dimensional subspaces of $W$. To see this, one has to show that if $(\zeta_n,v_n)\to(\zeta,v)$, $(\zeta'_n,v'_n)\to(\zeta',v')$ and $(\zeta_n,v_n)\sim(\zeta'_n,v'_n)$, then $(\zeta,v)\sim(\zeta',v')$. This is easy and left to the reader. Since different $C$-lines through 0 meet only at 0, we also see that $w_n\to w\ne0$ implies that $Cw_n\to Cw$.

It is clear that every $C$-line through 0 can be written in the form $C(\one,v)$ (uniquely), or $C(0,v)$ (non-uniquely). The $C$-lines of the first type form a dense open set (dense because $C(0,v)=\lim C(\frac1n,v)=\lim_{n\to\infty} C(\one,nv)$ as $n\to\infty$).

We call $GL(W,C)$ the group of $\R$-linear transformations of $W$ preserving $C$-lines. We also set $K=GL(W,C)\cap O(W)$. Since $GL(W,C)$ and $K$ are closed subgroups of $GL(W,\R)$, they are Lie groups.

\proclaim{Proposition 4.1} If $(\ph,\psi)\in M$, the map $\ph\times\psi:W\rightarrow W$ defined by $(\ph\times\psi)(\zeta,v)=\big(\ph(\zeta),\psi( v)\big)$ is in $K$. This correspondence identifies $M$ with the subgroup of $K$ whose elements fix the point $(\one,0)$. The subgroup $M_1$ of $M$ is then identified with the subgroup of $K$ whose elements fix every point in $C$.
\endproclaim

\proof If $(\ph,\psi)\in M$, then $\ph(\one)=\one$, so that $\ph\times\psi$ fixes $(\one,0)$. Let  $\ell=C(\zeta,v)$ be a $C$-line through 0. If $\zeta=0$, then $\ell=\{(0,\eta v):\eta\in C\}$, and 
$$
(\ph\times\psi)(0,\eta v)=\big(0,\psi(\eta v)\big)=\big(0,\ph(\eta)\psi(v)\big)\in C(\ph\times\psi)(0,v)\ .
$$

If $\zeta\ne0$, then $\ell=C(\one,\zeta\inv v)$, so we may as well assume that $\zeta=\one$. Then $\ell=\{(\eta,\eta v):\eta\in C\}$ and
$$
(\ph\times\psi)(\eta,\eta v)=\big(\ph(\eta),\ph(\eta)\psi(v)\big)\in C\big(\one,\psi(v)\big)= C(\ph\times\psi)(\one,v)\ .
$$

Conversely, assume that $k\in K$ fixes $(\one,0)$. Then $kC$ is the $C$-line containing $(\one,0)$, i.e. $C$. By orthogonality, $kV=V$. Denote by $k_1,k_2$ the restrictions of $k$ to $C$ and $V$ respectively. Consider now an element $(\eta,\eta v)$ in a $C$-line $C(\one,v)$ not contained in $V$. By linearity, $k(\eta,\eta v)=\big(k_1(\eta),k_2(\eta v)\big)$. Moreover, $k(\eta,\eta v)\in Ck(\one,v)=C\big(\one,k_2(v)\big)$. Then necessarily $k_2(\eta v)=k_1(\eta)k_2(v)$, i.e. $(k_1,k_2)\in M$.

The last part of the statement is now obvious.
\endproof

We will construct elements of $K$ that are not in $M$, i.e. that do not fix $(\one,0)$. But first we prove a more general statement, which will also be useful later.

\proclaim{Proposition 4.2} Let $\{v_1,\dots,v_n\}$ be an orthonormal $C$-basis of $V$. Let $w_0=(\one,0)$, $w_j=(0,v_j)$ ($1\le j\le n$) and $(a_{jk})_{0\le j,k\le n}$ be a real matrix. Then the map $A:W\rightarrow W$ defined by
$$
A(\sum_{j=0}^n\zeta_jw_j)=\sum_{j,k=0}^n a_{jk}\zeta_kw_j
$$
is in $GL(W,C)$.
\endproclaim

Though the expression $\zeta w$ does not make sense in general, the abuse of notation in the statement above does not cause ambiguity. It must be understood that $\zeta_0w_0=(\zeta_0,0)$ and $\zeta_jw_j=(0,\zeta_jv_j)$ for $j\ge1$.

\proof It is clear that $A$ is well defined and $\R$-linear, so that it suffices to prove the last statement for $C$-lines through 0. By continuity, we can also restrict ourselves to $C$-lines $\ell=C(\one,v)$. Now, for $v=\sum_{j=1}^n\zeta_jv_j$,
$$
C(\one,v)=\big\{(\la,\sum_{j=1}^n\la\zeta_jv_j):\la\in C\big\}=\big\{\sum_{j=0}^n(\la\cdot_{v_j}\zeta_j)w_j:\la\in C\big\}
$$
(where $\zeta_0=\one$ and $\la\cdot_{v_0}\zeta_0$ stands for $\la$). So,
$$
A(\ell)=\big\{\sum_{j,k=0}^na_{jk}(\la\cdot_{v_k}\zeta_k)w_j:\la\in C\big\}\ .
$$

This is clearly a $C$-line if $C$ is associative, because $\cdot _{v_k}$ is independent of $k$. 

If $C$ is non-associative, then $n=1$ and we must verify that
$$
\align
(a_{00}\la&+a_{01}\la\cdot_{v_1}\zeta_1)\inv(a_{10}\la+a_{11}\la\cdot_{v_1}\zeta_1)\\
&=\big[(a_{00}\one+a_{01}\zeta_1)\inv\cdot_{v_1}\la\inv\big]\cdot_{(a_{10}\la+a_{11}\la\cdot_{v_1}\zeta_1)v_1}\big[\la\cdot_{v_1}(a_{10}\one+a_{11}\zeta_1)\big]
\endalign
$$
does not depend on $\la$. For fixed $\la$, every element of $C$ appearing in this expression belongs to the division subalgebra of $(C,\cdot_{v_1})$ generated by $\zeta_1$ and $\la$. By Proposition 1.2 and Corollary 1.3, this subalgebra is associative and $\cdot_{(a_{10}\la+a_{11}\la\cdot_{v_1}\zeta_1)v_1}=\cdot_{v_1}$, so that $\la$ cancels out.
\endproof

\proclaim{Proposition 4.3} For a fixed unit vector $v_0\in V$, decompose $v\in V$ as
$v=\eta v_0+v'$, with $v'\in (Cv_0)^\perp$, and set, for $\theta\in\T$,
$$
\sigma_{v_0,\theta}(\zeta,\eta v_0+v')=\big(\cos\theta \zeta-\sin\theta\eta,
(\cos\theta\eta+\sin\theta\zeta)v_0+v'\big)\ .
$$

Then $\sigma_{v_0,\theta}\in K$.
\endproclaim

\proof It is trivial to verify that $\sigma_{v_0,\theta}\in O(W)$. It maps $C$-lines into $C$-lines because $v_0$ can be completed to an o.n. basis of $V$ and then Proposition 4.2 applies.
\endproof

\proclaim{Corollary 4.4} $K$ acts transitively on the unit sphere $S_W$ (and hence on $C$-lines through~$0$). 
\endproclaim

\proof Given a point $(\zeta,v)\in S_W$, write $(\zeta,v)$ as $(c\zeta',sv')$ with $\zeta'\in S_C$, $v'\in S_V$, $c=\cos\theta$, $s=\sin\theta$ for an appropriate $\theta\in\T$ (if $v=0$, we choose $v'$ arbitrarily).
Then $k=\sigma_{{\zeta'}\inv v',\theta-\frac\pi2}\circ \sigma_{v',\frac\pi2}$ maps $(\one,0)$ into $(\zeta,v)$.
\endproof

\vskip.2cm

We call  {\it linear $C$-subspace} of $W$ an $\R$-linear subspace $E$ such that $Cw\subset E$ whenever $w\ne0$ is in~$E$. 

We also call {\it affine $C$-subspace} of $W$ a translate $E'=E+w$ of a linear $C$-subspace $E$. This is equivalent to saying that $E'$ is a set that contains the whole affine $C$-line $w_1+C(w_2-w_1)$ connecting any pair of points $w_1,w_2\in E'$.

For an $\R$-linear subspace of $V$, to be a sub-$C$-module of $V$ is the same as being a linear $C$-subspace of $W$.

\proclaim{Lemma 4.5} The orthogonal of a linear $C$-subspace is a linear $C$-subspace.
\endproclaim

\proof Let $E$ be a linear $C$-subspace. Modulo the action of $K$, we can assume that $C\subset E$. Then $E=C\oplus V_0$, with $V_0$ a sub-module of $V$. Then the orthogonal of $E$ in $W$ is the same as the orthogonal of $V_0$ in $V$, which is also a sub-module. 
\endproof

\proclaim{Proposition 4.6} The linear span in $W$ of $C$-lines through $0$ is a linear $C$-subspace. Conversely, any linear $C$-subspace is the orthogonal sum of $C$-lines through~$0$.
\endproclaim

\proof The first part of the statement follows by induction, if we prove that the sum of a proper linear $C$-subspace $E$ and a $C$-line $\ell$ through 0 is a linear $C$-subspace.

Modulo the action of $K$, we can assume that $E^\perp$ contains $C$, hence that $E\subset V$. If $\ell$ is also contained in $V$, then $E+\ell$ is a sub-module of $V$.

If $\ell\not\subset V$, then $\ell=C(\one,v_0)$ for some $v_0\in V$ and 
$$
E+\ell=\big\{(\eta,\eta v_0+v):\eta\in C,v\in S\big\}\ .
$$

We need to prove that, if $\eta\ne0$ and $v\in E$, then $C(\eta,\eta v_0+v)\subset E+\ell$. An element of 
 $C(\eta,\eta v_0+v)$ has the form $(\zeta, \zeta v_0+\zeta\eta\inv v)$, which is in $E+\ell$, being $\zeta\eta\inv v\in E$.

For the converse, we can assume that $C\subset E$. Then $E=C\oplus V_0$, with $V_0$ a sub-$C$-module of $V$, and it can be decomposed into the orthogonal sum of irreducible ones.
\endproof

\vskip.2cm

Proposition 4.6 allows us to introduce the notions of {\it orthonormal $C$-basis} and of {\it $C$-dimension} of a linear $C$-subspace of $W$. We can then extend the scope of Corollary 4.4.

\proclaim{Corollary 4.7} The group $K$ acts transitively on the manifold $\Xi_W$ of ordered orthonormal $C$-bases of $W$ and on the manifold of linear $C$-subspaces of any fixed dimension.
\endproclaim

\proof It follows from Proposition 2.2 that $M\subset K$ acts transitively on ordered orthonormal $C$-bases of $W$ whose first element is $(\one,0)$.  It is then sufficient to prove that any other basis $\{w_0,\dots,w_n\}$ can be mapped into a basis of this type.  But this amounts to saying that there is $k\in K$ such that $kw_0=(\one,0)$, and this follows from Corollary 4.4.
\endproof

In the remainder of this section we describe the orbits of $M$ on the unit sphere $S_W$ of $W$ and analyze the structure of $K$ in more detail.

We assume that $C'\ne0$ and $V\ne0$, the degenerate cases being trivial. We fix an o.n. $C$-basis $\{v_1,\dots,v_n\}$ of $V$ and an element $z\in C'$ with $|z|=1$.

\proclaim{Lemma 4.8} Every $M$-orbit in $S_W$ meets the subspace $\R\one+\R z+\R v_1$ in the points $a\one\pm bz\pm cv_1$ for unique numbers $a\in\R$, $b,c\ge0$ with $a^2+b^2+c^2=1$.
\endproclaim

\proof $M$ fixes $\one$, hence $\RE\zeta$ remains constant on any orbit $M\cdot(\zeta,v)$. The rest of the statement follows from Corollary 2.3.
\endproof

We write, for $t\in\R$, $\eta_t=(\cos t)\one+(\sin t)z$, and define $m_t$ by
$$
m_t(\zeta,v)=(\eta_t\zeta \eta_t\inv,\eta_tv)\ .
$$

This is a one-parameter subgroup of $M$ ($m_t\in M_2$ if $\dim C'>1$, and $m_t\in M_1$ if $\dim C'=1$). We write $\sigma_1=\sigma_{v_1,\frac\pi2}$ and define $T$ as the one-parameter subgroup of $K$ consisting of the elements 
$$
\rho_t=\sigma_1\inv\circ m_t\circ \sigma_1\ .
$$

We also denote by $L$ the subgroup of $K$ preserving $C$ (hence also preserving $V$). Then $M\subset L$, and also $T\subset L$, since $\sigma_1$ only interchanges $C$ with $Cv_1$ and preserves every other $Cv_j$.

For any $v\in V$ and $t\in\R$, we have $\rho_t(\one,v)=(z_t,v')$ for some $v'\in V$. Hence the orbits of $L$ on $S_W$ all meet the plane $\R\one+\R v_1$ (in points $\pm a\one\pm cv_1$, for unique $a,c\ge0$ with $a^2+c^2=1$).

We write $T'$ for the subgroup $\{\sigma_{v_1,\theta}\}_{\theta\in\T}$ of $K$.

\proclaim{Proposition 4.9} We have $K=LT'L$ and $L=MTM$.
\endproclaim

The proof is trivial from the structure of the $L$- and $M$-orbits in $S_W$.

\proclaim{Corollary 4.10} The following hold:
\roster
\item"(i)" $M$ and $T$ generate $K$;
\item"(ii)" $MkM=Mk\inv M$ for every $k\in K$;
\item"(iii)" the projection of the action of $L$ onto $C$ is $SO(C)$ or $O(C)$\footnote{The projection of the action of $L$ onto $C$ is $SO(C)$ if $d=4,8$ and $O(C)$ if $d=1,2$.}.
\endroster
\endproclaim

\proof (i) is obvious and (ii) follows from Lemma 4.8. As to (iii), the description of $M$ given in Section 2 shows that $M_2$ acts on $C$ as $SO(C')$ or $O(C')$. On the other hand, $T$ acts on $C$ as a torus transversal to $SO(C')$. It is well known, and easy, that $SO(d-1)$ is a maximal proper subgroup of $SO(d)$, and the conclusion follows.
\endproof

\vskip.2cm
\subhead Remark 1 \endsubhead 
For future reference, we observe that, given $w,w'\in W$ and writing $\pi_{Cw}$ for the orthogonal projection on $Cw$, the ratio 
$|\pi_{Cw}w''|\big/|w''|$ is constant for all $w''\in Cw'$. This constant is then the cosine of the angle of $Cw$ and $Cw'$. The proof is immediate when $Cw=C$. The general case follows by Corollary 4.7.

\vskip.2cm
\subhead Remark 2 \endsubhead 
Corollary 4.10(ii) implies that $(K,M)$ is a Gelfand pair. This and (iii) are the key facts in the discussion of harmonic analysis on $K/M$ in the exceptional case developed in \cite{Ta}. Here we have proved them without explicit use of the octonions.

\vskip.5cm

\head 5. The construction of $CPW$ \endhead

\vskip.3cm

The compact rank-one symmetric spaces will be defined as appropriate compactifications of the vector spaces $W$ that we have associated to $J^2C$-modules. In the associative case, these compactifications can be described (as they are in the literature) as spaces of lines:  according to Theorem 3.1, we identify $C$ with an associative division algebra $\F$, $V$ with $\F^n$ and $W$ with $\F^{n+1}$. The symmetric space is then $\F P^{n+1}=(\F^{n+2}\setminus\{0\})/\F^*$. This construction makes use of the fact that we can ``add dimensions'' to $V$ without destroying the $J^2C$-module structure, something that cannot be adapted to the non-associative case.

In order to have a unified description, including both the associative and the non-associative case, we construct the compactification by ``gluing'' to $W$ a space $W_\infty$ consisting of ``points at infinity''. 

\vskip.3cm

We set $W_\infty=(W\setminus\{0\})/\sim$, where $\sim$ is the equivalence relation introduced  at the beginning of Section 4. It is convenient to think of $W_\infty$ as the set of $C$-lines through 0 in $W$. We denote by $\pi:W\setminus\{0\}\rightarrow W_\infty$ the quotient map, and we endow $W_\infty$ with the quotient topology. The element $\pi(w)$ will be denoted by $[w]$.

We then set $CPW=W\cup W_\infty$, and define a topology on it by assigning neighborhood bases at the various points as follows:
\roster
\item"(i)" to a point $w\in W$, we assign its Euclidean neighborhoods in $W$;
\item"(ii)" to a point $[w]\in W_\infty$, we assign, for each neighborhood $U_{[w]}^\infty$ of $[w]$ in $W_\infty$ and each $R>0$, the neighborhood
$$
U_{[w]}^R=U_{[w]}^\infty\cup\big\{w'\in\pi\inv (U_{[w]}^\infty):|w'|>R\big\}\ .
$$
\endroster

We leave the reader the verification that these neighborhood systems actually define a second-countable topology. It is obvious that $W$ is open and dense in $CPW$. 

As we have seen, a $C$-line not contained in $V$ contains a unique element $(\one,v)$. We denote the corresponding element of $W_\infty$ as $[\one,v]$ and $W_\infty^1$ the set of such elements. The complementary subset $W_\infty^0$ of $W_\infty$ is the set of $C$-lines contained in $V$. We use the notation $[0,v]$ to denote the $C$-line containing $v\ne0$ (in this case, however, $[0,v]$ depends only on the equivalence class  of $v$ mod. $C^*$).  Clearly, $W_\infty^1$ is open and dense in $W_\infty$.

The following facts about convergence of sequences in $CPW$ are easily derived from the definition:
\roster
\item"(i)" $W\ni w_n\to [w]\in W_\infty$ if and only if $|w_n|\to\infty$ and $[w_n]\to[w]$ in $W_\infty$;
\item"(ii)" $(\zeta_n,v_n)\to[\one,v]\in W_\infty^1$ if and only if $|\zeta_n|\to\infty$ and $\zeta_n\inv v_n\to v$;
\item"(iii)" $(\zeta_n,v_n)\to[0,v]\in W_\infty^0$ if and only if $|v_n|\to\infty$, $|\zeta_n|=o(|v_n|)$ and $\exists\ \la_n(\to0)$ in $C$ such that $\la_n v_n\to v$;
\item"(iv)" $W_\infty^1\ni[\one,v_n]\to[0,v]\in W_\infty^0$ if and only if $|v_n|\to\infty$ and $\exists\ \la_n(\to0)$ such that $\la_n v_n\to v$.
\endroster

We state a couple of facts that will be used later on.

\proclaim{Proposition 5.1} Let $f\in GL(W,C)$. The induced map $f_\infty:W_\infty\rightarrow W_\infty$ such that $f_\infty\circ\pi=\pi\circ f$ is a homeomorphism and $\bar f=f\cup f_\infty$ is a homeomorphism of $CPW$.

Translations in $W$ extend to homeomorphisms of $CPW$ which are the identity on $W_\infty$.
\endproclaim

\proof That $f_\infty$ is a homeomorphism follows from standard facts about quotient topologies. It is then sufficient to prove that $f(w_n)\to f_\infty([w])$ if $W\ni w_n\to[w]\in W_\infty$. 
Using (i) above, $|f(w_n)|\to\infty$ because $f$ is invertible, and 
$$
\big[f(w_n)\big]=f_\infty\big([w_n]\big)\to f_\infty([w])\ .
$$

Take now a translation $\tau(w)=w+w_0$ with $w_0=(\zeta_0,v_0)$, and assume that $(\eta_n,u_n)\to [w]\in W_\infty$. We claim that also $(\eta_n+\zeta_0,u_n+v_0)\to [w]$. 

If $[w]=[\one,v]\in W_\infty^1$, by (ii) we are assuming that $|\eta_n|\to\infty$ and $v_n=\eta_n\inv u_n\to v$. Then
$|\eta_n+\zeta_0|\to\infty$ and
$$
\align
\lim_{n\to\infty}(\eta_n+\zeta_0)\inv(u_n+v_0)&=\lim_{n\to\infty}(\eta_n+\zeta_0)\inv u_n\\
&=\lim_{n\to\infty}(\eta_n+\zeta_0)\inv \eta_nv_n\\
&=v-\lim_{n\to\infty}(\eta_n+\zeta_0)\inv \zeta_0v_n\\
&=v\ .
\endalign
$$

If $[w]=[0,v]\in W_\infty^0$, by (iii) we are assuming that $|u_n|\to\infty$, $|\eta_n|=o\big(|u_n|\big)$ and $\la_n u_n\to v$ for some sequence $\la_n\to0$. Then $|u_n+v_0|\to\infty$, $|\eta_n+\zeta_0|\sim|\eta_n|=o\big(|u_n+v_0|\big)$, and $\la_n(u_n+v_0)\to v$.
\endproof

From now on we will use the same symbol for a $C$-linear or $C$-affine map of $W$ and for its continuous extension to $CPW$.

\proclaim{Proposition 5.2} If $E$ is a linear $C$-subspace, the closure in $CPW$ of the affine $C$-subspace $E'=w_0+E$ is $(w_0+E)\cup\pi(E)$.
\endproclaim

\proof By Proposition 5.1, we can assume that $w_0=0$. Since $\pi(E)=\pi(E\cap S_W)$, $\pi(E)$ is closed in $W_\infty$. By (i), if $w_n\to[w]$ and $w_n\in E$, then $[w_n]\to[w]$, hence $[w]\in\pi(E)$. Conversely, if $[w]\in\pi(E)$, then $nw\to[w]$. 
\endproof

In order to introduce a differentiable structure on $CPW$ compatible with this topology, we first show that certain maps on $CPW$ are homeomorphisms. These maps will then be used to define the coordinate patches.

\proclaim{Lemma 5.3} Let $\ph_0:CPW\rightarrow CPW$ be defined as
$$
\matrix
\ph_0(\zeta,v)=(\zeta\inv,\zeta\inv v)\hfill&\text{if }\zeta\ne0\\
\ph_0(0,v)=[\one,v]\hfill&\\
\ph_0[\one,v]=(0,v)\hfill&\\
\ph_0[0,v]=[0,v]\ .\hfill&
\endmatrix\tag5.1
$$

Then $\ph_0$ is a homeomorphism.
\endproclaim

\proof Since $\ph_0\inv=\ph_0$, it is sufficient to prove that $\ph_0$ is continuous. The following facts must be verified:
\roster
\item $(\zeta_n,v_n)\to (0,v)$ (with $\zeta_n\ne0$) if and only if $(\zeta_n\inv,\zeta_n\inv v_n)\to[\one,v]$;
\item if $(\zeta_n,v_n)\to[0,v]$ (with $\zeta_n\ne0$), then $(\zeta_n\inv,\zeta_n\inv v_n)\to[0,v]$;
\item $(0,v_n)\to[0,v]$ if and only if $[\one,v_n]\to[0,v]$.
\endroster

These can be easily verified on the basis of (ii)-(iv).
\endproof

Fix now an orthonormal $C$-basis $\{v_1,\dots,v_n\}$ of $V$, and define, for $j=1,\dots,n$ and $v'\perp Cv_j$,
$$
\psi_j(\zeta,\eta v_j+v')
=(\eta,\zeta v_j+v')\ .
$$

Then $\psi_j$ is linear and involutive on $W$ and maps $C$-lines into $C$-lines. Therefore it extends to an involutive homeomorphism, also denoted by $\psi_j$, of $CPW$.
If $\ph_j=\ph_0\circ \psi_j\circ \ph_0$, then also $\ph_j$ is an involutive homeomorphism of $CPW$. Explicitly, with $v'\perp Cv_j$,
$$
\matrix
\ph_j(\zeta,\eta v_j+v')&=(\eta\inv\cdot_{v_j}\zeta,\eta\inv v_j+\eta\inv v')\hfill&\text{if }\eta\ne0\hfill\\
\ph_j(\zeta,v')\hfill&=\left\{\matrix
\big[\one,\zeta\inv(v_j+v')\big]\hfill&\text{if }\zeta\ne0\ ,\hfill\\[0,v_j+v'] \hfill&\text{if }\zeta=0\
,
\hfill\endmatrix\right.\ \hfill
\\
\ph_j[\one,\eta
v_j+v']\hfill&=\left\{\matrix
(\eta\inv,
\eta\inv v')\hfill&\text{if }\eta\ne0\ ,\hfill\\  [\one,v'] \hfill&\text{if }\eta=0\ ,\hfill
\endmatrix\right.\hfill\\ 
\ph_j[0,\eta
v_j+v']\hfill&=\left\{\matrix
(0,\eta\inv v')\hfill&\text{if }\eta\ne0\ ,\hfill\\[0,v'] \hfill&\text{if }\eta=0\ .\hfill
\endmatrix\right.\hfill
\endmatrix\tag5.2
$$

For notational convenience, we set $\ph_{n+1}=\id_{CPW}$.

\proclaim{Proposition 5.4} The charts $\big(\ph_j(W),\ph_j\big)$, with $0\le j\le n+1$, define a differentiable structure on $CPW$.
\endproclaim

(Since $\ph_j=\ph_j\inv$, the coordinate maps take values in $W$.)

\proof
We show that $W,\ph_0(W),\dots \ph_n(W)$ form an open covering of $CPW$. In fact, the points of $W_\infty$ that are not contained in $\ph_0(W)$ are those in $\pi(V)$, whereas, for $j\ge1$, the points of $W_\infty$ not contained in $\ph_j(W)$ are those in $\pi\big((Cv_j)^\perp\big)$. Since
$$
V\cap (Cv_1)^\perp\cap\cdots\cap (Cv_n)^\perp=0\ ,
$$
we conclude that
$$
\bigcup_{j=0}^n \ph_j(W)\supset W_\infty\ .
$$

We see directly from (5.1) and (5.2) that each component of $\ph_j$, $0\le j\le n$, is a rational function. The same is then true also for the transition maps $\ph_j\circ\ph_k$. In particular they are smooth (and analytic) on $W\cap(\ph_j\circ\ph_k)\inv W$. 
\endproof

When restricted to $W_\infty$, this differentiable structure coincides with the quotient structure of $(W\setminus\{0\})/\sim$. Restricting the quotient map $\pi$ to the unit sphere $S_W$, we obtain  the {\it Hopf fibration}
$$
\pi:S_W\rightarrow W_\infty\ ,
$$
with fiber $S_C$. By Proposition 5.1, $K$ acts by diffeomorphisms of $W_\infty$ and, by Corollary 4.7, this action is transitive. Let $L$ be the subgroup of $K$ introduced in Section 4, whose elements preserve $C$ as a set. Identifying $W_\infty$ with $K/L$ and the unit sphere $S_W$ with $K/M$, it follows that $L$ is the structure group of the bundle, and the fiber $S_C$ is diffeomorphic to $L/M$.

\vskip.3cm

It may help the reader to see explicitly how $CPW$ identifies with $\F P^{n+1}$, when $\F$ is an associative division algebra, $C=\F$ and $V=\F^n$.  
Let $\Pi:(\F^{n+2}\setminus\{0\})\rightarrow\F P^{n+1}$ be the quotient map. We think of $W$ inside $CPW$ as $\Pi\big(\F^{n+1}\times\{1\}\big)$ and of $W_\infty$ as $\Pi\big(\F^{n+1}\times\{0\}\big)$. Then each $\ph_j$ is the projective image of the map interchanging the $j$-th component of $(q_0,q_1,\dots,q_n,q_{n+1})\in\F^{n+2}$ with the $(n+1)$-th component.

\vskip.5cm

\head 6. $CPW$ as a compact symmetric space \endhead

\vskip.3cm

We first introduce a metric on $W$ and then prove that it can be extended to all of $CPW$. Guided by the non-compact case (cf. (A.1) in the Appendix), we set
$$
\lan X,Y\ran_{w+}=\left\{\matrix \displaystyle\frac{\lan X,Y\ran}{(1+|w|^2)^2}\hfill&\text{ if }X,Y\in
Cw\ ,\hfill\\ 
\displaystyle\frac{\lan X,Y\ran}{1+|w|^2}\hfill&\text{ if }X,Y\in
(Cw)^\perp\ ,\hfill\\
\displaystyle0\hfill& \text{ if } X\in Cw, Y\in (Cw)^\perp\ ,\hfill
\endmatrix\right.\tag6.1
$$
if $w\ne0$ and, for $w=0$, 
$$
\lan X,Y\ran_{0+}=\lan X,Y\ran\qquad \text{ for every }X,Y\ .\tag6.2
$$

Notice that (6.1) and (6.2) imply that for arbitrary $X\in W$,
$$
|X|_{w+}=c_{|w|,\ph}|X|\ ,\tag6.3
$$
where $\ph$ is the angle of the $C$-lines $Cw$ and $CX$ (see Remark 1 at the end of Section~4) and
$$
c_{|w|,\ph}=\Big(\frac{\cos^2\ph}{(1+|w|^2)^2}+\frac{\sin^2\ph}{1+|w|^2}\Big)^\half\ .\tag6.4
$$

\proclaim{Proposition 6.1} The geodesic $\gamma$ going through the origin with tangent vector $w\in S_W$ is  $\gamma(t)=(\tan t) w$, $|t|<\frac\pi2$.
The group of isometries of $W$ fixing the origin is $K$.
\endproclaim

\proof The elements of $K$ are isometries of $W$, as a direct consequence of their property of mapping $C$-lines into $C$-lines. We can then assume that $w=(\one,0)$. Then $\gamma$ must be invariant under $M$. By Kostant's double transitivity (Corollary 2.3), this implies that $\gamma\subset \R\one$, i.e. $\gamma(t)=r(t)w$. We impose now that 
$$
\|\Dot\gamma(t)\|_{\gamma(t)+}=\frac{r'(t)}{1+r^2(t)}=1\ .
$$

This gives $r(t)=\tan t$.

Finally, let $f$ be an isometry of $W$ with $f(0)=0$, and let $k=Df(0)\in O(W)$. If $w\in S_W$, the geodesic $\gamma(t)=(\tan t) w$ is mapped into $f\big(\gamma(t)\big)=(\tan t) k(w)=k\big(\gamma(t)\big)$. Therefore $f=k$. Then (6.1) implies that $k$ must map $C$-lines into $C$-lines.
\endproof

We will see that the transition maps $\ph_j\circ \ph_k$ are isometric on $W\cap(\ph_j\circ\ph_k)\inv(W)$, and this will allow us to extend the metric to $CPW$, defining it locally as the pull-back of the metric in $W$ via the coordinate maps $\ph_j$. Since $\psi_j\in K$ for $1\le j\le n$, the crucial fact to be proved is that $\ph_0$ is an isometry on $W\setminus V$.

This will be the first consequence of the next lemma. 

\proclaim{Lemma 6.2} For $\theta\in \T$ and $\zeta\ne(\cot \theta)\one$, define
$$
b_\theta(\zeta,v)=\big((\cos\theta\one-\sin\theta\zeta)\inv(\sin\theta\one+\cos\theta \zeta),
(\cos\theta\one-\sin\theta\zeta)\inv v\big)\ ,\tag6.5
$$
where the product in the first component is in the sense of Clifford multiplication or, equivalently, with respect to any product $\cdot_{v_0}$.
Then $b_\theta$ is an isometry.
\endproclaim

Before giving the proof, we comment that the $b_\theta$ must be regarded as close analogues of the maps 
$\sigma_{v_0,\theta}$ introduced in Proposition 4.3. We can see this in the associative case, where $CPW=\F P^{n+1}$. In terms of the homogeneous coordinates $(q_0,\dots,q_{n+1})$ introduced at the end of the previous section, and with $v_0$ the $j$-th basis element ($1\le j\le n$), $\sigma_{v_0,\theta}$ introduces a rotation by $\theta$ in the pair of coordinates $(q_0,q_j)$. Now $b_\theta$ introduces the same rotation in the coordinates $(q_0,q_{n+1})$.

We also remark the formal analogy between (6.5) and (A.3) in the Appendix.

\proof Fix $w=(\zeta,v)$ and a tangent vector $(\eta,u)$ at $w$. Set $b_\theta(w)=w_\theta=(\zeta_\theta,v_\theta)$.  We compute 
$$
(b_\theta)_{*,w}(\eta,u)=\frac d{d\eps}_{|_{\eps=0}} \Big(\big(c\one-s(\zeta+\eps\eta)\big)\inv\big(s\one+c (\zeta+\eps\eta)\big),
\big(c\one-s(\zeta+\eps\eta)\big)\inv (v+\eps u)\Big)\ ,
$$
 where $\eps\in\R$, $c=\cos\theta$, $s=\sin\theta$, and prove that $\|(b_\theta)_{*,w}(\eta,u)\|^2_{w_\theta+}$ does not depend on $\theta$ as long as $\zeta\ne(\cot\theta)\one$. It is convenient to restrict ourselves to $\zeta\not\in\R\one$, $v\ne0$, the general case following by continuity.
 
 We take the first-order expansion of each component in $\eps$. The expansion of the first component takes place in the subalgebra of $(C,\cdot_{v_0})$ (for some $v_0$) generated by $\zeta$ and $\eta$, which is associative. By standard computations, we find that
 $$
 \aligned
 &(b_\theta)_{*,w}(\eta,u)\overset{\text{def}}\to=(\eta_\theta,u_\theta)\\
 &=\big((c\one-s\zeta)\inv\eta(c\one-s\zeta)\inv,(c\one-s\zeta)\inv u+s(c\one-s\zeta)\inv\eta(c\one-s\zeta)\inv v\big)\ .
 \endaligned\tag6.6
 $$

Notice that, according to Remark 2 in Section 1, the expression defining $\eta_\theta$ is independent of the choice of $v_0$.

The subspace $C+Cw$ and its orthogonal complement $(C+Cw)^\perp$ (with respect to the product inner product on $W$) are clearly invariant under $ (b_\theta)_{*,w}$. Notice that, since $w_\theta$ is in $C+Cw$, then $C+Cw$ is also the direct sum of $Cw_\theta$ and its orthogonal complement. By (6.1), $C+Cw$ and $(C+Cw)^\perp$ are also orthogonal with respect to the inner product $\lan\ ,\ \ran_{w_\theta+}$. It is then sufficient to prove that $\|(\eta_\theta,u_\theta)\|^2_{w_\theta+}$ does not depend on $\theta$  when $(\eta,u)$ belongs to either of the two subspaces.

Take $(0,u)\in (C+Cw)^\perp$, i.e. with $u\perp Cv$ in $V$. Then $\eta_\theta=0$ and $u_\theta=(-s\zeta+c\one)\inv u$. Then
$$
\|(0,u_\theta)\|_{w_\theta+}^2=\frac{|u_\theta|^2}{1+|w_\theta|^2}=\frac{|u|^2}{(1+|w_\theta|^2)|-s\zeta+c\one|^2}\ .
$$

Observing that
$$
1+|w_\theta|^2=\frac{1+|w|^2}{|c\one-s\zeta|^2}\ ,\tag6.7
$$
the conclusion follows.

Take now $\xi=(\eta,\la\zeta\inv v)\in C+Cw$. Observe that $\xi\in Cw$ if and only if $\la=\eta$ and $\xi\perp Cw$ if and only if $\la=-\frac{|\zeta|^2}{|v|^2}\eta$. It follows that the orthogonal decomposition of $\xi$, for general $\eta,\la$, into a term in $Cw$ and one orthogonal to it is
$$
(\eta,\la\zeta\inv v)=(\al,\al\zeta\inv v)+\Big(\beta,-\frac{|\zeta|^2}{|v|^2}\beta\zeta\inv v\Big)
$$
with
$$
\align 
\al&=\frac{|\zeta|^2}{|w|^2}\eta+\frac{|v|^2}{|w|^2}\la\ ,\\
\beta&=\frac{|v|^2}{|w|^2}(\eta-\la)\ .
\endalign
$$

Therefore, skipping the straightforward computations,
$$
\aligned
\big\|(\eta,\la\zeta\inv v)\big\|^2_{w+}&=\frac{|\al|^2\Big(1+\frac{|v|^2}{|\zeta|^2}\Big)}{(1+|w|^2)^2} +
\frac{|\beta|^2\Big(1+\frac{|\zeta|^2}{|v|^2}\Big)}{1+|w|^2} \\
&=\frac1{(1+|w|^2)^2}\big(|\eta|^2+|\la\zeta\inv v|^2+|\eta-\la|^2|v|^2\big)\ .
\endaligned\tag6.8
$$

Obviously, the same formula holds with $w,v,\zeta$ replaced by $w_\theta,v_\theta,\zeta_\theta$.
We then apply (6.8) for a given $\theta$, with $\eta=\eta_\theta$ as in (6.6) and $\la=\la_\theta$ such that $(\eta_\theta,\la_\theta\zeta_\theta\inv v_\theta)=(b_\theta)_{*,w}(\eta,\la\zeta\inv v)$, i.e.,
$$
\la_\theta\zeta_\theta\inv v_\theta=(c\one-s\zeta)\inv\big(\la\zeta\inv +s\eta(c\one-s\zeta)\inv\big) v\ .\tag6.9
$$

By (6.7),
$$
\frac{|\eta_\theta|^2}{(1+|w_\theta|^2)^2}=\frac{|\eta|^2}{(1+|w|^2)^2}\ ,
$$
which does not depend on $\theta$. The next term is
$$
\frac{|\la_\theta\zeta_\theta\inv v_\theta|^2}{(1+|w_\theta|^2)^2}=\frac{|c\one-s\zeta|^2\big|\big(\la\zeta\inv +s\eta(c\one-s\zeta)\inv\big) v\big|^2}{(1+|w|^2)^2}\ .
$$

For the last term, we use the identity
$$
|\eta_\theta-\la_\theta|^2|v_\theta|^2=|\zeta_\theta|^2|\eta_\theta-\la_\theta|^2|\zeta_\theta\inv v_\theta|^2
=|\zeta_\theta|^2|(\eta_\theta-\la_\theta)\zeta_\theta\inv v_\theta|^2\ .
$$

Simple calculations and the identity $ (s\one+c\zeta)\inv-s\one=c(c\one-s\zeta)(s\one+c\zeta)\inv$ give that
$$
(\eta_\theta-\la_\theta)\zeta_\theta\inv v_\theta=(c\one-s\zeta)\inv\big(c\eta(s\one+c\zeta)\inv-\la\zeta\inv\big)v\ .
$$

Hence
$$
\frac{|\eta_\theta-\la_\theta|^2|v_\theta|^2}{(1+|w_\theta|^2)^2}=\frac{|s\one+c\zeta|^2\big|\big(c\eta(s\one+c\zeta)\inv -\la\zeta\inv\big)v \big|^2}{(1+|w|^2)^2}\ .
$$

We must then prove that
$$
|c\one-s\zeta|^2\big|\big(\la\zeta\inv +s\eta(c\one-s\zeta)\inv\big) v\big|^2+|s\one+c\zeta|^2\big|\big(c\eta(s\one+c\zeta)\inv -\la\zeta\inv\big)v \big|^2
$$
does not depend on $\theta$. Expanding the right-hand factor in each summand, this quantity equals
$$
\align
\big(|&c\one-s\zeta|^2+|s\one+c\zeta|^2\big)|\la\zeta\inv v|^2+(s^2+c^2)|\eta|^2|v|^2\\
&\qquad-2|c\one-s\zeta|^2s\big\lan\la\zeta\inv v,\eta(c\one-s\zeta)\inv v\big\ran-2|s\one+c\zeta|^2c\big\lan\la\zeta\inv ,v\eta(s\one+c\zeta)\inv v\big\ran\\
&=(1+|\zeta|^2)|\la\zeta\inv v|^2+|\eta|^2|v|^2\\
&\qquad+2s\big\lan\la\zeta\inv v,\eta(c\one-s\bar\zeta) v\big\ran-2c\big\lan\la\zeta\inv ,v\eta(s\one+c\bar\zeta) v\big\ran\\
&=(1+|\zeta|^2)|\la\zeta\inv v|^2+|\eta|^2|v|^2-2\lan\la\zeta\inv v,\eta\bar\zeta v\ran\ ,
\endalign
$$
a quantity that does not depend on $\theta$.
\endproof

\proclaim{Corollary 6.3} The transition maps $\ph_j\circ\ph_k$ are isometries on $W\cap(\ph_j\circ\ph_k)\inv W$.
\endproclaim

\proof Since $\ph_0=(-\id)\circ b_{\frac\pi2}$, and $-\id\in K$, $\ph_0$ is an isometry on $W\setminus V$. Since $\psi_j\in K$ for $1\le j\le n$, then $\ph_j$ is an isometry on $W\setminus(Cv_j)^\perp$. Compositions are then isometries on appropriate open dense subsets of $W$ and, by continuity, they remain isometric on $W\cap(\ph_j\circ\ph_k)\inv W)$.
\endproof

We can then extend the metric to $CPW$ by imposing that the $\ph_j$, $0\le j\le n$, are isometries of the whole space. 
Observe that
\roster\item"(i)" by Proposition 6.1, $W_\infty$ is the geodesic sphere centered at the origin of $W$ and radius $\frac\pi2$;
\item"(ii)" the $b_\theta$ extend uniquely to isometries of $CPW$; for $\theta=0,\pi$, $b_\theta\in K$, so that its extension is obvious; for $\theta\ne0,\pi$ we have
$$
\aligned
b_\theta\big((\cot\theta)\one,v\big)&=[\one,(\sin\theta)v]\ ,\\
b_\theta\big([\one,v]\big)&=\big(-(\cot\theta)\one,-(\sin\theta)\inv v\big)\ ,\\
 b_\theta\big([0,v]\big)&=[0,v]\ ;
 \endaligned\tag6.10
$$
\item"(iii)" $\gamma(\theta)=b_\theta(0)$ is the geodesic through 0 with tangent vector $(\one,0)$ there, and its length is $\pi$.
\endroster

\proclaim{Theorem 6.4} $CPW$ is a compact symmetric space of rank one. 
\endproclaim

\proof We show that $CPW$ is homogeneous by proving that the origin can be mapped to any other point by isometries. If $w\in W$, take $\theta=\arctan|w|$, so that $b_\theta(0)=(|w|,0)$. By Corollary 4.4, there is now $k\in K$ such that $k(|w|,0)=w$. Take now $p=\pi(Cw)\in W_\infty$, and let $k\in K$ be such that $kC=Cw$. Then the extension $\bar k$ of $k$ to $CPW$ maps $[\one,0]$ into $p$. 
By (5.1), $\bar k\circ\ph_0$ maps 0 into $p$.

The geodesic symmetry around 0 is $-\id$, which is in $K$, hence $CPW$ is symmetric. Since the action of $K$ on the unit sphere in the tangent space at 0 is transitive, the rank of $CPW$ is one. 
\endproof

 Let $U$ be the isometry group of $CPW$ and $B=\{b_\theta:\theta\in\T\}$. 
 
 \proclaim{Corollary 6.5} $U=KBK$.  The centralizer of $B$ in $K$ is $M$.
 \endproclaim
 
 \proof The proof of Theorem 6.4 shows that  any isometry of $CPW$ can be decomposed as $k_1b_\theta k_2$ with $k_1,k_2\in K$.
 Let $k\in K$ be such that $b_\theta k=kb_\theta$ for every $\theta$. Then
 $$
 \big((\tan\theta)\one,0\big)=b_\theta k(0,0)=kb_\theta(0,0)=k\big((\tan\theta)\one,0\big)\ ,
 $$
i.e. $k$ fixes the points in $(\one,0)$. By Proposition 4.1, $k\in M$. Conversely, take $m\in M$. Then $m=(\ph,\psi)$, where $\ph,\psi$ satisfy (2.1). Given $(\zeta,v)\in W$, let $(\zeta',v')=\big(\ph(\zeta),\psi(v)\big)$. Set $c=\cos\theta, s=\sin\theta$. Hence
$$
b_\theta(\zeta',v')=\big((c\one-s\zeta')\inv(s\one+c\zeta'),(c\one-s\zeta')\inv v'\big)
$$

It follows from (2.1) that 
$$
\ph(\eta\cdot_v\la)=\ph(\eta)\cdot_{\psi(v)}\ph(\la)\ .
$$

By Remark 1 in Section 1, if $\eta$ and $\la$ are rational expression in $\zeta$, the value of the product does not depend on $v$. We can then say that
$$
(c\one-s\zeta')\inv(s\one+c\zeta')=\ph\big((c\one-s\zeta)\inv(s\one+c\zeta)\big)\ ,
$$
and
$$
(c\one-s\zeta')\inv v'=\ph\big((c\one-s\zeta)\inv\big)\psi(v)=\psi\big((c\one-s\zeta)\inv v\big)\ ,
$$
i.e., $b_\theta(\zeta',v')=mb_\theta(\zeta,v)$.
\endproof

\vskip.5cm

\head 7. Some applications \endhead

\vskip.3cm

The purpose of this section is to use our setup to prove a few (known) fundamental facts about compact symmetric spaces of rank one. In the usual approach these things are easy to prove for spheres and classical projective spaces, but quite difficult for the octonionic case (cf. \cite{Be}, Ch. 3).

Throughout this section we consider a $J^2C$-module $(C,V)$ and the associated $W=C\oplus V$, together with $CPW$. If $V$ has a $C$-basis of $n$ elements, we say that $\dim_CW=n+1$.

\proclaim{Theorem 7.1} If $d=\dim_\R C$ and $m=\dim_CW$, we have
$$
\vol(CPW)=\frac{\Gamma\big(\frac d2\big)}{\Gamma\big((m+1)\frac d2\big)}\pi^{m\frac d2}\ .
$$
\endproclaim

\proof Since $W$ is open dense in $CPW$, by (6.1) we have
$$
\vol(CPW)=\int_W(1+|w|^2)^{-2d-(m-1)d}\,dw\ ,
$$
where $dw$ is Lebesgue measure. Polar coordinates and the substitution $x=|w|^2$ reduce this to
$$
\frac12|S^{md-1}|\int_0^\infty (1+x)^{-(m+1)\frac d2}x^{m\frac d2-1}\,dx\ ,
$$
which is a classical Beta-integral.
\endproof

\proclaim{Theorem 7.2} Let $E$ be a plane in the tangent space to $CPW$ at $(0,0)$ (identified with $W$ as usual), and suppose that $X,Y\in W$ form an orthonormal basis of $E$. Then the sectional curvature $\sigma(E)$ of $E$ is
$$
\sigma(E)=1+3|\pi_{CX}Y|^2=1+3\cos^2\ph\ ,\tag7.1
$$
where $\pi_{CX}$ denotes orthogonal projection onto $CX$ and $\ph$ is the angle of $CX$ and $CY$.
\endproclaim

\proof We use the classical formula 
$$
\sigma(E)=\lim_{r\to0}\frac3\pi \frac{2\pi r-L_r}{r^3}\ ,\tag7.2
$$
where $L_r$ is the arc length of the exponential of a circle of radius $r$ in $E$. By Proposition 6.1, this curve can be written as 
$$
\gamma_r(\theta)=\tan r\big((\cos\theta)X+(\sin\theta)Y\big)\ ,\qquad (0\le\theta\le2\pi)\ .
$$

Since $\gamma_r(\theta)$ and $\overset.\to{\gamma_r}(\theta)=\tan r\big(-(\sin\theta)X+(\cos\theta)Y\big)$ are images of $\gamma_r(0)=(\tan r)X$ and $\overset.\to{\gamma_r}(0)=(\tan r)Y$ respectively under a real rotation, the angle $\ph$ of $C\gamma(\theta)$ and $C\overset.\to{\gamma_r}(\theta)$ is independent of $\theta$ and $\cos\ph=|\pi_{CX}Y|$ (cf. Remark 1 at the end of Section 4). Using the abbreviation $p=|\pi_{CX}Y|$ and (6.3), (6.4), it follows that
$$
L_r=2\pi\frac{\tan r}{1+\tan^2r}\big(1+(1-p^2)\tan^2r\big)^\half\sim 2\pi\Big(r-\big(\frac16+\frac{p^2}2\big)r^3\Big)\ ,
$$
and (7.2) gives the result.
\endproof

\subhead Remark  \endsubhead
The $C$-lines of $W$ have now the following geometric characterization: two linearly independent elements $X,Y$ in $W$ are in the same $C$-line through 0 if and only if $\sigma(E)=4$ for the plane $E$ spanned by them. Similarly, $CX\perp CY$ if and only if $\sigma(E)=1$. 

\vskip.3cm

Next, we describe the totally geodesic submanifolds of $CPW$. We fix an o.n. $C$-basis $\{v_1,\dots,v_n\}$ of $V$ and take a number $n_0$, $0\le n_0\le n$. If $n_0=0$, let $C_0$ be any real subspace of $C$. If $n_0\ge1$, let $C_0$ be any division subalgebra of $(C,\cdot_{v_1})$. In either case we set
$$
W_0=C_0\oplus \sum_{j=1}^{n_0}C_0v_j=C_0\oplus V_0\ .\tag7.3
$$

We note that $(C_0,V_0)$ is a $J^2C$-module. If $n>1$, this follows from Proposition~1.4 and Corollary 1.5. If $n=1$, it is trivial.

\proclaim{Theorem 7.3} The closure of $W_0$ in $CPW$ is totally geodesic. All totally geodesic submanifolds of $CPW$ arise as $U$-images of these.
\endproclaim

\proof It is well known and obvious that the fixed point set of an isometry in a Riemannian manifold is totally geodesic. The linear map which is the identity on $W_1=C\oplus\sum_{j=1}^{n_0}Cv_j$ and minus the identity on $\sum_{j=n_0+1}^n Cv_j$ is in $M$, so it is an isometry. Therefore $W_1$ is totally geodesic. It is then sufficient to prove that $W_0$ is totally geodesic in $W_1$.

If $n_0=0$, there is not much to say. In this case $W_1=C$. If $C_0$ is any subspace of $C$, the reflection with respect to $C_0$ is an isometry, and this implies that $C_0$ is totally geodesic. The closure of $C$ in $CPW$ is a $d$-dimensional sphere, and the closure of $C_0$ is a lower-dimensional geodesic sphere.

Assume therefore that $n_0\ge1$.
Given a division subalgebra $C_0$ of $C$, there is a chain of division subalgebras between $C_0$ and $C$, each of index 2 in the next. So it is enough to consider the case where $C_0$ has index 2 in $C$. In this case we have the orthogonal direct sum $C=C_0\oplus C_0z$, where $z$ is an arbitrary unit element orthogonal to $C_0$ (the product being understood as $\cdot_{v_1}$). 

Define $\al:C\rightarrow C$ by $\al(\zeta_1+\zeta_2z)=\zeta_1-\zeta_2z$, for $\zeta_1,\zeta_2\in C_0$. Observe that left multiplication by a unit element $\eta\in C$ is orthogonal on $C$, which preserves the decomposition $C_0\oplus C_0z$ if $\eta\in C_0$, and interchanges the two summands if $\eta\in C_0z$. Using this, it follows that $\alpha$ is an automorphism of $C$. 

We extend $\al$ to $W_1$ by $\al(\zeta, \sum \zeta_jv_j)=\big(\al(\zeta), \sum \al(\zeta_j)v_j\big)$. This is an orthogonal transformation and it preserves $C$-lines through the origin. In fact, if $\zeta\inv\zeta_jv_j={\zeta'}\inv\zeta'_jv_j $ for every $j$, then also $\al(\zeta)\inv\al(\zeta_j)v_j=\al({\zeta'})\inv\al(\zeta'_j)v_j $, in the associative case obviously, and in the non-associative case because there is only one $j$.

Since $W_0$ is the fixed point set of $\al$, this proves that $W_0$ is totally geodesic in $W_1$, hence in $W$. Because the metric induced on $W_0$ from $W$ is the same metric constructed starting from the $J^2C$-module $(C_0,V_0)$, the embedding of $W_0$ in $CPW$ extends to an isometric embedding of $C_0PW_0$, whose image is the closure of $W_0$. This proves the first statement.

To prove the converse, let $N\subset CPW$ be totally geodesic. $N$ is a symmetric space and clearly it has rank one. Translating by an element of $U$, we may assume that $0\in N$. The tangent space to $N$ at 0 is a subspace $W_0$ of $W$. Because geodesics through the origin are straight lines, $N$ is the closure of $W_0$ in $CPW$. 

By Theorem 8.4, whose proof is independent of the present theorem, there is a $J^2C$-module $(C_0,V_0)$ such that $N$ is isometric (up to a constant factor) to $C_0PW_0$. Because $C_0$ has constant curvature equal to 4, it follows from the Remark above that $\dim C_0\le\dim C$.  Acting by $K$, we can assume that $C_0\subset C$. Similar considerations based on the same Remark imply that $V_0\subset V$. Along the same lines, if $\{v_1,\dots, v_{n_0}\}$ is an orthonormal $C_0$-basis of $V_0$, the subspaces $C_0v_j$ are contained in different, hence orthogonal, $C$-lines in $V$. Therefore $\{v_1,\dots, v_{n_0}\}$ can be completed to an orthonormal $C$-basis $\{v_1,\dots, v_n\}$ of $V$, and $W_0$ has the form~(7.3).
\endproof

We call {\it projective $C$-line} in $CPW$ a maximal totally geodesic submanifold of constant curvature equal to 4. If no sectional curvature in $CPW$ is equal to 4 (i.e. if $\dim C=1$), we call projective $C$-lines the geodesics.

\proclaim{Corollary 7.4} $U$ acts transitively on projective $C$-lines. The projective $C$-lines intersecting $W$ are the closures in $CPW$ of affine $C$-lines in $W$.  The projective $C$-lines contained in $W_\infty$ are the images $\pi(E)$ of linear $C$-subspaces $E$ of $W$ of $C$-dimension two.
\endproclaim

\proof The first statement follows immediately from Theorem 7.3. 

Observe now that, for $0\ne v\in V$, the affine $C$-line $(\tan\theta,0)+Cv$ in $W$ is the image of $Cv$ under the isometry $b_\theta$ in (6.5). Therefore it is a geodesic in $W$, and the same is true for its closure in $CPW$. By Corollary 6.5, this gives the description of all projective $C$-lines intersecting $W$.

One projective $C$-line entirely contained in $W_\infty$ is obtained as the image under $\ph_0$ of the projective $C$-line $Cv\cup\{[0,v]\}$ for a fixed $v\ne0$. This gives $\{[\one,\zeta v]:\zeta\in C\}\cup\{[0,v]\}$, equal to $\pi(C\oplus Cv)$. If we act on this line by an element $u=k_1b_\theta k_2$ of $U$, the image will remain inside $W_\infty$ if and only if $b_\theta\pi\big(k_2(C\oplus Cv)\big)\subset W_\infty$. By (6.10), if $E$ is a $C$-subspace of dimension 2 and $b_\theta$ maps $\pi(E)$ into $W_\infty$, then necessarily $\theta\in\{0,\pi\}$, so that $u\in K$.

It follows that all the projective $C$-lines in $W_\infty$ are the images of $C\oplus Cv$ under the action of $K$, and the conclusion follows from Corollary 4.7.
\endproof

It is clear that any two distinct points of $CPW$ determine a unique projective $C$-line.
We also note that the projective $C$-lines are the ``Helgason spheres'' which, by \cite{Hel2}, are present in any compact symmetric space.

\vskip.3cm

Generalizing the notion of projective $C$-line, we call {\it projective $C$-subspace} a subset $E$ of $CPW$ such that for any pair of distinct points of $E$, the whole projective $C$-line joining them is contained in $E$. If a projective $C$-subspace $E$ intersects $W$, it is easy to see that $E\cap W$ is an affine $C$-subspace of $W$. Since we can always put ourselves in this situation applying an element of $U$, it follows that $E$ is a submanifold of real dimension equal to a multiple $kd$ of $d$. We then say that $E$ has $C$-dimension equal to $k$.
By Theorem 7.3, the projective $C$-subspaces are the $U$-images of $W_0\cup\pi(W_0)$, with $W_0$ as in (7.3) and $C_0=C$. We also have the following analogue of Corollary~7.4.

\proclaim{Corollary 7.5} $U$ acts transitively on $k$-dimensional projective $C$-subspaces. The $k$-dimensional projective $C$-subspaces intersecting $W$ are the closures in $CPW$ of affine $C$-subspaces in $W$.  The $k$-dimensional projective $C$-subspaces contained in $W_\infty$ are the images $\pi(E)$ of linear $C$-subspaces $E$ of $W$ of $C$-dimension $k+1$.
\endproclaim

\vskip.3cm

A number of further geometric facts follow very easily. For instance, the Jacobi fields along a geodesic can be determined by reduction to the geodesic $\gamma(t)=(\tan t,0)$. For $z\in C'$ and $v\in V$, let $z(t)$ and $v(t)$ be their parallel displacements along $\gamma$. Because $C$ is totally geodesic in $W$, $z(t)$ is the same as parallel displacement in $C$ of $z$ along $\gamma$. Since $C$ has constant curvature 4, it follows (cf. \cite{KN, vol. I, p. 71}) that $(\sin 2t)z(t)$ and $(\cos 2t)z(t)$ are Jacobi fields. Similarly, $v(t)$ is  the same as parallel displacement in $\R\one\oplus \R v$ of $v$ along $\gamma$. Therefore $(\sin t)v(t)$ and $(\cos t)v(t)$ are Jacobi fields. By a dimension count, all Jacobi fields along $\gamma$ are linear combination of these and of $\overset.\to\gamma(t)$ and $t\overset.\to\gamma(t)$.

\vskip.3cm

Another fact concerns cut locus and conjugate points of a given point $p\in CPW$. We can assume that $p=0$. As we have already observed, two geodesics, $\exp_0(tX)$ and $\exp_0(tY)$, meet only at 0 if $X$ and $Y$ are not in the same $C$-line, and at the two points 0 and $[X]=[Y]$ if $X$ and $Y$ belong to the same $C$-line. It is then clear that the cut locus of 0 is $W_\infty$, or, in other words, the set of points $q$ such that $d(0,q)=\frac\pi2$. In the latter formulation, the statement is true for any $p$ in place of 0.

It also follows that the conjugate locus of any point $p$ is $\{p\}$ when $d=\dim C=1$ (the real projective space), and coincides with the cut locus otherwise. In the latter case the multiplicity of the first conjugate point along any geodesic is $d-1$.

\vskip.3cm

Finally we note that $CPW$ has a natural cell decomposition 
$$
CPW={\Cal C}^{d(n+1)}\cup{\Cal C}^{dn}\cup\cdots \cup{\Cal C}^0\ ,
$$
where ${\Cal C}^k$ denotes a cell of dimension $k$. In fact, we can take ${\Cal C}^{d(n+1)}=W$. Its complement $W_\infty$ is isometric, under either of the maps $\ph_j$ in (5.2), with the closure of $C\oplus (Cv_j)^\perp=W_j$. This closure is $CPW_j$, with one less $C$-dimension. So the statement follows by induction.

\vskip.5cm

\head 8. Every compact rank-one symmetric space is a $CPW$ \endhead

\vskip.3cm

The reader who is willing to accept the classification list of symmetric spaces existing in the literature, e.g. \cite{Be, Hel1, W}, can compare that list with the list in Section 3 of $J^2C$-modules and the corresponding $CPW$'s, and convince himself that we have obtained all compact rank-one symmetric spaces. However, in this section we will give a classification-free proof of this fact, clarifying at the same time the duality relations between compact and non-compact spaces.

\proclaim{Lemma 8.1} $CPW$ is simply connected if and only if $d=\dim C>1$.
\endproclaim

\proof If $d=1$, then $CPW$ is $\R P^{n+1}$ if $n=\dim V>1$ and $S^1$ if $V=0$. Therefore $CPW$ is not simply connected. 

Suppose now that $d>1$, and observe that if $\dim W>1$ and $W_\infty$ is simply connected, so is $CPW$. This depends on the fact that $W_\infty$ is a deformation retract of $CPW$ minus a point in $W$. 

If $V$ is non-trivial and $v_1$ is a unit vector in $V$, set $V'=(Cv_1)^\perp$. 
The map $\ph_1$ in (5.2) establishes a diffeomorphism between $W_\infty$ and $CPW'$ (cf. end of Section 7), with $W'=C\oplus V'$. By induction, matters are reduced to $V=0$. In this case $CPW$ is the sphere $S^d$, which is simply connected.
\endproof

It is well known (cf. \cite{W}) that the duality between compact and non-compact orthogonal semisimple Lie algebras induces a (bijective) duality between simply connected compact symmetric spaces and symmetric spaces of the non-compact type. Duality respects the rank, and the balls $B_W$ described in the Appendix, with the metric defined in (A.1), give us models of all the rank-one symmetric spaces of non-compact type. In fact, we can restrict ourselves to those $W=C\oplus V$ with $d=\dim C>1$.

\proclaim{Proposition 8.2} Let $W=C\oplus V$ with $d>1$. The  simply connected compact dual of the unit ball $B_W$ is $CPW$.
\endproclaim

\proof Assume that $d>1$. The non-compact dual of $CPW$ must be a ball $B_{W^\sharp}$ for some $W^\sharp=C^\sharp\oplus V^\sharp$ with the same dimension as $W$. Take the origin in $W$ as base points in $CPW$ and the origin in $W^\sharp$ as base point in $B_{W^\sharp}$.

By \cite{W, Cor. 8.4.3}, the two tangent spaces can be identified in such a way that sectional curvatures of corresponding planes have opposite values, possibly up to a constant factor. It follows from Theorem 7.2 and the comments following Lemma A.1 in the Appendix that $C$ and $C^\sharp$ must have the same dimension. Then $W^\sharp=W$.
\endproof

It follows that the $CPW$ with $d>1$ are all the simply connected compact symmetric spaces of rank one. Setting aside the trivial one-dimensional case ($C=\R$, $V=0$), every other compact symmetric spaces of rank one must have a $CPW$ with $d>1$ as its simply connected covering. We are so led to discuss the existence of locally isometric symmetric quotients of the $CPW$.

\proclaim{Lemma 8.3} Suppose that $CPW$ is a non-trivial covering of a symmetric space $X$. Then $V=0$ and $CPW$ is a two-fold covering of $X$.
\endproclaim

\proof Let $\pi:CPW\rightarrow X$ be the covering map. Since $X$ has rank one, its geodesics are circles and they all have the same length.
Let $E$ be the set of points $p\in CPW$, different from 0, that are mapped to $\pi(0)$. Then $E$ is non-empty and finite. Take $p_0\in E$ of minimal distance from 0 and let $\gamma$ be a full geodesic circle in $CPW$ going through 0 and $p$. Then $\pi(\gamma)$ is a geodesic circle in $X$ and its length is $\delta=d(0,p_0)$, strictly smaller than the length of $\gamma$. It follows that every geodesic circle in $CPW$ containing 0 is mapped by $\pi$ onto a geodesic circle of length $\delta$, hence it must contain a point in $E$ at distance $\delta$ from 0. 

Since $E$ is finite, we surely have two distinct geodesics through 0 intersecting at a point $p\in E$. By Proposition 5.2, geodesics through zero have the form $\R w\cup[w]$ with $w\in W$. Therefore two distinct geodesics can have a common point different from 0 only if they lie in the same $C$-line $Cw$, and in this case the common point is $[w]$. This implies that $\delta=\frac\pi2$ and hence $[w]\in E$ for every $w\in W$. So $W_\infty$ is finite. But this is only possible if $W=C$, and then $E=W_\infty$ consists of one single point.
\endproof

\proclaim{Theorem 8.4} The $CPW$ are all the distinct compact symmetric spaces of rank one, including the circle $S^1$.
\endproclaim

\proof It remains to verify that if $d=1$ and $W=\R\oplus V$ with $V\ne0$, then $CPW$ is isometric, up to a factor, to the quotient of $CP\tilde W$ modulo the antipodal map, with $\tilde C=W$ and $\tilde V=0$. To see this, consider the map
$\pi$ from $\tilde W$ to $W$ given by
$$
\pi(w)=\frac{2w}{1-|w|^2}\ ,\qquad \big(|w|\ne1\big)\ .
$$

This is an isometry up to a factor 2, and it extends continuously to the unit sphere and to the point at infinity, identifying the pairs of antipodal points, $w$ and $-\frac w{|w|^2}$.
\endproof

\vskip.5cm

\head 9. The group $GL(W,C)$ \endhead

\vskip.3cm

We study $GL(C,W)$ in some detail and, to avoid trivialities, we assume that $V\ne0$.

$GL(W,C)$ is closed under adjoints. To see this, note that for any $\R$-subspace $W_1$, $w\perp gW_1$ if and only if $g^*w\perp W_1$. Hence $g^*\big((gW_1)^\perp\big)=W_1^\perp$. Now let $Cw$ be a $C$-line and let $W_1=g\inv\big((Cw)^\perp\big)$. By Lemma 4.5, $(Cw)^\perp$ is a $C$-subspace, and by Proposition 4.6 so is $W_1$. By Lemma 4.5 again, $g^*(Cw)$ is a $C$-line.

It follows that $GL(W,C)$ is a reductive Lie group, since it is fixed under the Cartan involution $g\mapsto {g^*}\inv$ of $GL(W,\R)$ (cf. \cite{W}). What follows is a useful characterization of the elements of $GL(W,C)$.

\proclaim {Theorem 9.1} Assume that $V\ne0$. The elements of $GL(W,C)$ are the homeomorphisms of $W$ onto itself fixing the origin and mapping parallel $C$-lines into parallel $C$-lines.
\endproclaim

\proof One implication is obvious. So let $g$ be a homeomorphism fixing 0 and mapping parallel $C$-lines into parallel $C$-lines. We need to prove that $g$ is $\R$-linear. Since $g$ is continuous, it suffices to prove that it is additive.

If $w,w'$ are points in $W$ belonging to different $C$-lines through 0, it is easily verified that
$w+w'$ is the only point in the intersection of $w+Cw'$ with $w'+Cw$. 

Let $w,w'$ be as above. By assumption, $g(Cw)=Cg(w)$ and $g(w'+Cw)=g(w')+Cg(w)$. Similarly, $g(Cw')=Cg(w')$ and $g(w+Cw')=g(w)+Cg(w')$. By injectivity, $Cg(w)$ and $Cg(w')$ are different $C$-lines.
Therefore $g(w+w')\in \big(g(w')+Cg(w)\big)\cap \big(g(w)+Cg(w')\big)$. By the previous remark, $g(w+w')=g(w)+g(w')$.

If $w,w'$ belong to the same $C$-line through 0, the same identity follows by conitnuity, using the existence of other $C$-lines.
\endproof

By Corollary 4.7, every $g\in GL(W,C)$ can be decomposed as $g=k\circ h$ with $k\in K$ and $h\in GL(W,C)$ such that $h(V)=V$.

\proclaim{Lemma 9.2} Assume that $h\in GL(W,C)$ maps $V$ into itself and let
$$
h=\pmatrix \al&0\\ \al(\cdot)v_0&\ph\endpmatrix\ ,\tag9.1
$$
be the matrix representing it relative to the decomposition $W=C\oplus V$.
Then 
\roster
\item"(i)" $\al:C\rightarrow C$ is a scalar multiple of an orthogonal transformation,
\item"(ii)"  $\ph:V\rightarrow V$ is invertible and satisfies $\ph(\zeta v)=\al(\zeta)\al(\one)\inv\ph(v)$,
\item"(iii)" $v_0\in V$.
\endroster

Conversely, every triple $(\al,\ph,v_0)$ satisfying {\rm (i), (ii), (iii)} defines, through $(9.1)$, an element of $GL(W,C)$ mapping $V$ into itself.
\endproclaim

\proof Any $h$ mapping $V$ into itself is represented by a matrix
$$
h=\pmatrix \al&0\\ \sigma&\ph\endpmatrix\ .
$$

The condition $h(C)$ equal to a $C$-line implies that 
$$
\al(\zeta)\inv\sigma(\zeta)=\al(\one)\inv\sigma(\one)= v_0\ ,
$$
i.e. $\sigma(\zeta)=\al(\zeta)v_0$. 

For $v\in V$, impose now that $h(0,\zeta v)\sim h(0,v)$. This gives 
$\ph(\zeta v)=\beta(\zeta,v)\ph(v)$, with $\beta(\zeta,v)\in C$ and $\beta(\one,v)=1$.

Imposing now that 
$$
\psi(\zeta,\zeta v)=\big(\al(\zeta),\al(\zeta)v_0+\beta(\zeta,v)\ph(v)\big)\sim\psi(\one,v)=\big(\al(\one),\al(\one)v_0+\ph(v)\big)\ ,
$$
we obtain that 
$$
v_0+\al(\zeta)\inv\beta(\zeta,v)\ph(v)=v_0+\al(\one)\inv\ph(v)\ ,
$$
i.e. $\beta(\zeta, v)=\al(\zeta)\al(\one)\inv$.

If $z\in C'$ and $|z|=1$, then $z^2=-\one$, and
$$
\ph(v)=-\ph(z^2v)=-\big(\al(z)\al(\one)\inv\big)^2\ph(v)\ .
$$

We then have $|\al(z)|=|\al(\one)|$ and $\al(z)\cdot_{\ph(v)}\al(\one)\inv\in C'$. In particular $\al(z)\perp\al(\one)$ for $z\in C'$. It follows that $|\al(\zeta)|=|\al(\one)|$ for every $\zeta\in C$.

To proof of the converse is easy and we leave it to the reader.
\endproof

An immediate consequence of this lemma is the following statement, which will be used in Section 10.

\proclaim{Corollary 9.3} Let $g\in GL(W,C)$. The restriction of $g$ to any $C$-line $Cw$ through the origin is a scalar multiple of an orthogonal map.
\endproclaim

\proof In order to prove this, decompose $g=kh$ with $k\in K$ and $h$ as in Lemma 9.2. It is then sufficient to assume that $g=h$. By density, we can also assume that $w=(\one,v)$. Take $w'=(\zeta,\zeta v)\in Cw$. Then $|w'|=|\zeta||w|$ and
$$
\align
|h(w')|^2&=|\al(\zeta)|^2+\big|\al(\zeta)v_0+\al(\zeta)\al(\one)\inv\ph(v)\big|^2\\
&=|\al(\zeta)|^2|\al(\one)|^{-2}\big(|\al(\one)|^2+\big|\al(\one)v_0+\ph(v)\big|^2\big)^2\\
&=|\zeta|^2|h(w)|^2\ .
\endalign
$$

Hence the quantity $|h(w')|/|w'|$ is constant on $Cw$.
\endproof

\vskip.3cm

We fix now an ordered orthonormal $C$-basis $\{w_0,\dots, w_n\}$ of $W$ with
$w_0=(\one,0)$ and the other $w_j=(0,u_j)$ in $V$. Then any element $w$ of $W$ can be written as
$$
w=\sum_{j=0}^n\zeta_jw_j\ ,
$$
with the same abuse of notation as in Proposition 4.2.
This basis induces a flag $\{W_j\}$ of
linear $C$-subspaces in
$W$, with 
$$
W_j=\span\{w_j,\dots,w_n\}\ .\tag9.2
$$ 

Let $P$ be the subgroup of $GL(W,C)$
consisting of the elements $h$ that preserve the flag, i.e. such that
$h(W_j)=W_j$ for every $j$.

We introduce three subgroups of $P$. The first group is $M_P=P\cap K=P\cap L$. By Corollary 4.7, $M_P$ acts transitively on the product of unit spheres in the various $Cw_j$. The second group $\A\cong \R^{n+1}$ acts as scalar multiplication by $t_j>0$ on each $Cw_j$. The third is the group $N$ of those $h\in P$ such that $h_{|_{Cw_j}}=\id\,(\text{\rm mod\,}W_{j+1})$ for each~$j$. 

It is easy to verify that the elements of $N$ are in one-to-one correspondence with the lower-triangular matrices with entries $\la_{ij}\in C$,
$$
\Lambda=\pmatrix 1&0&\dots&\dots&0\\ \lambda_{10}&1&0&&0\\ \lambda_{20}&\lambda_{21}&1&\ddots&\vdots\\ \vdots&\vdots&\ddots&\ddots&\vdots\\ \lambda_{n0}&\lambda_{n1}&\dots&\dots&1\endpmatrix\ ,\tag9.3
$$
in the sense that to each $\Lambda$ as above we associate $n\in N$ given by
$$
n\Big(\sum_{j=0}^n\zeta_jw_j)=\sum_{j=0}^n\Big(\sum_{k<j}\zeta_k\cdot_{v_j}\lambda_{jk}+\zeta_j\Big)w_j\ .
$$

The product $n_1n_2$ in $N$ corresponds to the matrix product $\trans(\trans\Lambda_2\trans\Lambda_1)$ (the double transposition is a consequence of the fact that we write the action of  $C$ on $V$ as a left action). Notice that, if $C$ is associative, multiplication in $C$ is unambiguously defined, whereas, if $C$ is non-associative, the matrices $\Lambda$ are $2\times2$ and their product only involves the sum of their $(1,0)$-entries.

We then have the following Langlands decomposition of $P$ and Iwasawa decomposition of $GL(W,C)$.

\proclaim{Proposition 9.4} $P$ is the semidirect product
$M_P \A N$. $M_P$ and $\A$
commute, and $M_P \A$ normalizes $N$. Moreover, $GL(W,C)=K\A N$.
\endproclaim

\proof It is clear that  $M_P$ and $\A$
commute, and that $M_P \A$ normalizes $N$.

Given $h\in P$, let $\al$, $\ph$, $v_0$ be as in Lemma 9.2. 

For $j\ge1$, $h(w_j)=\big(0,\ph(u_j)\big)\in W_j$. Invertibility of $\ph$ implies that the  $(Cu_j)$-component of $\ph(u_j)$ must be different from zero.
We can then set 
$$
\ph(u_j)=\eta_j(u_j+\lambda_{j+1,j}u_{j+1}+\cdots +\lambda_{nj}u_n)\qquad(\forall\,j\ge1)\ ,\tag9.4
$$
with $\eta_j\ne0$.
If
$$
\aligned
w'_0&=\big(\one,v_0)\ ,\\
w'_j&=(0,u_j+\lambda_{j+1,j}u_{j+1}+\cdots +\lambda_{nj}u_n)\qquad(\forall\,j\ge1)\ ,
\endaligned\tag9.5
$$
then $h$ maps each $C$-line $Cw_j$ onto $Cw'_j$. Let
$$
v_0=\lambda_{10}u_1+\cdots +\lambda_{n0}u_n\ ,\tag9.6
$$
and let  $\Lambda$ be the lower-triangular matrix with all 1's along the diagonal and the $\la_{ij}$  defined in (9.4) and (9.6). The corresponding element $n\in N$ also maps each $C$-line $Cw_j$ onto $Cw'_j$. Therefore $n\inv h$ maps $Cw_j$ into itself for every $j$.

Precisely, one verifies that
$$
\align
n\inv h(\zeta,0)&=\big(\al(\zeta),0\big)\ ,\\
n\inv h(0,\zeta u_j)&=\big(0,\al(\zeta)\al(1)\inv\eta_ju_j\big)\qquad(\forall\,j\ge1)\ .
\endalign
$$

It follows that $n\inv h$ is a scalar multiple of an orthogonal transformation on each $C$-line $Cw_j$, hence $n\inv h=ma\in M_P \A$. So $h=nma$, and, since $M_P \A$ normalizes $N$, $h=man'$ for some other $n'\in N$.

To show uniqueness, it suffices to observe that  $M_P\cap \A$
and $(M_P \A)\cap N$ are trivial.

We have already pointed out that any element of $GL(W,C)$ is the product of an element of $K$ and an element of $P$. Therefore $g\in GL(W,C)$ can be decomposed as $kan$, with $k\in K$, $a\in \A$, $n\in N$. If $kan=k'a'n'$, then $k\inv k'=an(a'n')\inv$ is the identity on each $Cw_j$, hence $k=k'$. Then $a=a'$ and $n=n'$, by the previous part of the proof.
\endproof

Finally, we easily obtain the Cartan decomposition of $GL(W,C)$.

\proclaim{Theorem 9.5} $GL(W,C)=K\A K$.
\endproclaim

\proof Given $g\in GL(W,C)$, $g^*g$ is $C$-linear, self-adjoint and positive. Its eigen\-spaces are $C$-subspaces, so there is an orthonormal $C$-basis $\{e_0,\dots,e_n\}$ of $W$ such that $g^*g=t_j^2\id$ on $Ce_j$ for every~$j$. Take $k\in K$ such that $ke_j=w_j$, with $w_j$ as above, and let $a\in \A$ be such that $a_{|_{Cw_j}}=t_j\id$ for every~$j$.  Then  $a\inv kg^*gk\inv a\inv=\id$, i.e. $gk\inv a\inv\in K$.
\endproof

\vskip.5cm

\head 10. The group of collineations \endhead

\vskip.3cm

We consider here a larger group of transformations of $CPW$ than the isometry group $U$, i.e. the group of homeomorphisms that preserve the class of projective $C$-lines. These transformations are called {\it collineations} in projective geometry. 
Clearly, this definition makes sense only if $CPW$ is not reduced to one single $C$-line, i.e. if $V\ne0$. As we will see in Corollary 10.4, the restriction of a collineation to a single $C$-line turns out to be a conformal map. For this reason, when $V=0$ the {\it conformal group} is the natural substitute of the collineation group. In the following we will assume that $V$ is non-trivial. 

The elements of $U$ are obviously collineations and such are translations by Corollary 7.4. We remark in the Appendix that the group of collineations also includes the isometry group of the non-compact dual $B_W$ of $CPW$. Collineations obviously map projective $C$-subspaces into projective $C$-subspaces of the same dimension.

Denote by $\G$ the collineation group. If $g$ is a collineation, there exist $u\in U$ such that $u(W_\infty)=g(W_\infty)$. Then there is a translation $\tau$ such that 
$$
g=u\tilde g \tau\ ,\tag10.1
$$ 
with $\tilde g$ a collineation fixing 0 and stabilizing $W_\infty$. It follows from Theorem 9.1 that $\tilde g$ is the continuous extension of an element of $GL(W,C)$.
More generally, we can say that the collineations stabilizing $W_\infty$ are the extensions of $C$-affine invertible transformations of $W$.

It follows from (10.1) that collineations map projective $C$-subspaces into projective $C$-subspaces of the same dimension. 

\vskip.3cm

Let $\{W_j\}$ be the flag  (9.2) of linear $C$-subspaces of $W$, and consider the chain of projective $C$-subspaces of $CPW$
$$
\pi(W_n)\subset \cdots\subset\pi(W_0)=W_\infty\ .
$$

Let $P_\G$ be the subgroup of $\G$ whose elements preserve each $\pi(W_j)$. Then (up to extension to $CPW$) $P_\G$ is the semidirect product of $P$, introduced in Section 5, with the translation group $\W\cong W$. Let $N_\G$ be the nilpotent group given by the semidirect product of $N\subset P$ with $\W$. We then have the following decomposition of $\G$.

\proclaim{Theorem 10.1} Any $g\in \G$ can be uniquely decomposed as $g=uan$, with $u\in U$, $a\in \A$ and $n\in N_\G$.
\endproclaim

\proof Let $g=u \tilde g \tau$ as in (10.1). By Proposition 9.4, $\tilde g=kan_0$ with $k\in K$, $a\in \A$, $n_0\in N$. Then $g=(uk)a(n_0\tau)\in U\A N_\G$.

If now $uan=u'a'n'$, with $u,u'\in U$, $a,a'\in \A$, $n,n'\in N_\G$, then $u\inv u'=an(a'n')\inv\in P_\G$. Writing $n=n_0\tau$, $n'=n'_0\tau'$, with $n_0,n'_0\in N$, $\tau,\tau'\in \W$, the fact that $an(a'n')\inv\in GL(W,C)$ implies that $\tau=\tau'$. Therefore $u\inv u'=an_0(a'n'_0)\inv\in P$. Hence $u\inv u'\in K\cap P$ and $u=u'$. Finally, $an_0=a'n'_0$ implies that $a=a'$, $n_0=n'_0$.
\endproof

We introduce now an involution on $\G$. 
As observed in Section 7, $W_\infty$ is the set of points at distance $\frac\pi2$ from the origin of $W$, also characterized as the cut locus of 0. By $U$-invariance, it follows that the cut locus $p^*$ of any point $p\in CPW$ is a projective $C$-hyperplane, and that every projective $C$-hyperplane is the cut locus of some point. The correspondence $p\mapsto p^*$ is a bijection between points and projective $C$-hyperplanes.

\proclaim{Lemma 10.2} The correspondence $p\mapsto p^*$ is such that
$$
\matrix
\hfill 0^*&=W_\infty\ ,\hfill &&\\
\hfill w^*\cap W&=(Cw)^\perp-\frac1{|w|^2}w\hfill &\qquad&(\forall\,w\in W\setminus\{0\})\ ,\hfill \\
\hfill [w]^*\cap W&=(Cw)^\perp\hfill &\qquad&(\forall\,[w]\in W_\infty)\ .\hfill 
\endmatrix\tag10.2
$$

To each collineation $g$ we can associate the map $\theta g:CPW\rightarrow CPW$ defined by
$$
\big(\theta g(p)\big)^*=g(p^*)\ .
$$

Then 
\roster
\item"(i)" $\theta g=(g^*)\inv$ if $g\in GL(W,C)$, 
\item"(ii)" for every $g$, $\theta g$ is a collineation, 
\item"(iii)" $\theta$ is an involutive automorphism of $\G$, 
\item"(iv)"  $\theta g=g$ if and only if $g\in U$.
\endroster
\endproclaim

\proof For $g\in U$, we trivially have that $g(p^*)=\big(g(p)\big)^*$ (i.e. $\theta g=g$). Modulo the action of $K$, it is then sufficient to identify $(t\one,0)^*$ for $t\in\R$ and $[\one,0]^*$ in order to obtain (10.2). Let $b_s$, $s\in\T$, be the maps defined in (6.5). Since $b_s(0)=(\tan s\,\one,0)$ for $s\ne\pm\frac\pi2$ and $b_{\pm\frac\pi2}(0)=[\one,0]$, the problem reduces to determining $b_s(W_\infty)$. By~(6.10),
$$
b_s(W_\infty^1)=\Big\{\Big(-(\cot s)\one, -\frac1{\sin s}v\Big):v\in V\Big\}=-\big((\cot s)\one\big)+V\ ,
$$
and $(t\one,0)^*$ is its closure. Taking $s=\frac\pi2$, this gives $[\one,0]^*$.

The map $\theta g$ is well defined for $g\in\G$. 
If $g\in U$, $g(p^*)=\big(g(p)\big)^*$, i.e.  $\theta g=g$. Moreover, for every $g,h\in\G$,
$$
\big((\theta g)(\theta h)(p)\big)^*=g\big((\theta h(p))^*\big)=gh(p^*)=\big(\theta(gh)(p)\big)^*\ ,
$$
i.e. $(\theta g)(\theta h)=\theta(gh)$.

To prove (i), fix $w\in W$, $w\ne0$, and take $k\in K$ such that $kw=(t\one,0)$ with $t>0$. Writing $g=g'k$, we have $\theta g(w)=(\theta g')kw=\theta g'(t\one,0)$.
Decompose $g'$ as $k'h$ with $k'\in K$ and $h$ as in (9.1). Then $\theta g'(t\one,0)=k'\theta h(t\one,0)$.

By (10.2),
$$
h\big((t\one,0)^*\big)=h\big(V-(t\inv\one,0)\big)=V-\big(t\inv\al(\one),0\big)\ ,
$$
i.e. $\theta h(t\one,0)=\big(t\al(\one)|\al(\one)|^{-2},0\big)=(h^*)\inv(t\one,0)$. Finally,
$$
\theta g(w)=\theta(k'hk)(w)=k'(\theta h)k(w)=k'(h^*)\inv k(w)=(g^*)\inv(w)\ .
$$

To prove (ii), by (10.1) it suffices to prove that, for any translation $\tau_{w_0}(w)=w+w_0$, $\theta\tau_{w_0}$ is a collineation. Since for $w\in W$, $w\ne0$,
$$
\tau_{w_0}(w^*)=(Cw)^\perp-\frac1{|w|^2}w+w_0=(Cw)^\perp-\frac1{|w|^2}w+P_{Cw}w_0\ ,
$$
with $P$ denoting orthogonal projection, $\theta\tau_{w_0}(w)$ is the element $w'\in Cw$ such that
$$
\frac1{|w'|^2}w'=\frac1{|w|^2}w-P_{Cw}w_0\ .
$$

Take $k\in K$ such that $kw_0=(t\one,0)$ with $t>0$. Then $\tau_{w_0}=k\inv\tau_{kw_0}k$ and $\theta\tau_{w_0}=k\inv(\theta\tau_{kw_0})kw$. We can then restrict ourselves to $w_0=(t\one,0)$.

We compute $P_{Cw}(\one,0)$ for $w=(\eta,u)$, assuming that $\eta\ne0$ (for $\eta=0$ the projection is 0, and $w'=w$). We impose that
$$
(\one-\la,-\la\eta\inv u)\perp (\al,\al\eta\inv u) \qquad (\forall\,\al\in C)
$$
and obtain that
$$
P_{Cw}(\one,0)=\frac1{|w|^2}\big(|\eta|^2,\bar\eta u\big)\ .
$$

With standard simplifications, this gives
$$
\theta\tau_{w_0}(\eta,u)=\big((\one-t\eta)\inv\eta,(\one-t\eta)\inv u\big)\ .\tag10.3
$$

To prove that this is a collineation, we observe that (10.3) admits a continuous bijective extension to $CPW$.
By continuity, we can limit ourselves to verify that the image of an affine $C$-line in $W$ that is not parallel to $V$ is contained in a projective $C$-line. Therefore we take 
$$
w=\om_\zeta=(\zeta_0+\zeta, v_0+\zeta v_1)\ ,
$$
with $\zeta$ varying in $C$ and $\zeta_0,v_0,v_1$ fixed. We show that the points 
$$
\theta\tau_{w_0}(\om_\zeta)-\theta\tau_{w_0}(\om_0)=(\la_\zeta, u_\zeta)\ ,
$$
with
$$
\align
\la_\zeta&=(\one-t(\zeta_0+\zeta))\inv(\zeta_0+\zeta)-(\one-t\zeta_0)\inv\zeta_0\\
u_\zeta&=(\one-t(\zeta_0+\zeta))\inv (v_0+\zeta v_1)-(\one-t\zeta_0)\inv v_0
\endalign
$$
lie on the same $C$-line through 0. A straightforward computation shows that $\la_\zeta\inv u_\zeta$ is independent of $\zeta$.

Concerning (iii), we have already observed in the course of this proof that $\theta$ is multiplicative. By (10.1), (i) and (ii), $\theta$ maps $\G$ into itself. The identity $\theta^2g=g$ holds for $g\in U$ obviously and for $g\in GL(W,C)$ by (i). Computing $\theta^2\tau_{w_0}$ from (10.3), we conclude that it also holds for translations. 

Finally, assume that $\theta g=g$. Modulo $U$, we can assume that $g(0)=0$. Therefore, $g(W_\infty)=\big(\theta g(0)^*\big)=W_\infty$. Therefore $g\in GL(W,C)$ and $g=(g^*)\inv$, i.e.~$g\in K$. 
\endproof

\proclaim{Theorem 10.3} $\G$ is a simple Lie group G with restricted root system of type $A_{n+1}$ and the differential of $\theta$ (also denoted by $\theta$) is a Cartan involution on its Lie algebra. The corresponding Cartan decomposition of the Lie algebra $\g$ of $\G$ is $\g=\u+\p$, with $\u$ the Lie algebra of $U$. The Lie algebra $\a$ of $\A$ is maximal abelian in $\p$ and the decomposition in Theorem 10.1 is the associated Iwasawa decomposition.
\endproclaim

\proof We sketch two proofs of the fact that $\G$ is a Lie group. The first one makes use of the theorem of Montgomery-Zippin-Gleason. Since $\G$ with the topology of uniform convergence on $CPW$ is clearly a topological transformation group, it suffices to prove that $\G$ is locally Euclidean.

The stabilizer of in $\G$ of $W_\infty$ as a set is the semidirect product $GL(W,C)\times_s\W$. Hence the stabilizer of 0 is its $\theta$-image $B=GL(W,C)\times_s\theta \W$. Since $\theta$ is easily seen to be a homeomorphism of $\G$, $B$  is locally Euclidean with the topology induced from $\G$. Since $\W\cdot0=W$, it follows that $\W B$ is a neighborhood of the identity in $\G$. If $U_\W$ and $U_B$ are small neighborhoods of the identity in $\W$ and $B$ respectively, then $U_\W U_B$ is a neighborhood of the identity in $\G$ and is homeomorphic to $U_\W\times U_B$, showing that $\G$ is locally Euclidean.

The second proof consists in introducing local coordinates near the identity coming from the decomposition $\W\cdot GL(W,C)\cdot\theta \W$ of an open dense subset of $\G$. The map $g\mapsto \theta g\inv$ is clearly smooth on this set, and one needs to verify smoothness of the group operations. 
The key  point is then to verify, with the aid of the formulas defining $\theta$, that the map $(w_0,w)\mapsto \theta\tau_{w_0}(w)$ of $W\times W$ into $W$ is rational, and therefore smooth where meaningful. From this one gets that $\theta$ and the inversion are both smooth on a neighborhood of the identity. With appropriate use of left translations, one finishes the proof.

The existence of the Cartan involution $\theta$ implies that $\g$ is reductive. In order to prove that it is semisimple, we must show that its center is trivial. The previous part of the proof shows that
$$
\g=\gl+\w+\theta\w\ ,\tag10.4
$$
where $\gl$ the Lie algebra of $GL(W,C)$ and $\w\sim W$ that of $\W$. 

Let $Z=X+t_w+\theta t_{w'}$ be a central element. For every $Y\in\gl$, we must have
$$
\align
0=[Z,Y]&=[X,Y]+[t_w,Y]+\theta [t_{w'},\theta Y]\\
&=[X,Y]+[t_w,Y]-\theta [t_{w'},Y^*]\ .
\endalign
$$

Since the decomposition (10.4) is respected, each term on the right-hand side must be zero.  In particular, $[Y,t_w]=0$ for every $Y$. Since $[Y,t_w]=t_{Yw}$, this implies that $w$ is fixed by all elements of $GL(W,C)$ in the connected component of the identity. Therefore $w=0$. Similarly, $w'=0$ so that $Z=X$. But then $[X,t_{w''}]=0$ for every $w''$, and the action of $\gl$ on $W$ is effective. Therefore $X=0$.

By (10.4), $\p=\p\cap\gl+(\id-\theta)\w$ and, by Proposition 9.4, $\p\cap\gl=\a+(\id-\theta)\n$, where $\n$ is the Lie algebra of $N$. 

The argument preceding Proposition 9.4 shows that, given $a=(a_0,\dots,a_n)\in\A$, $ana\inv=(a_i/a_j)n$ when $n\in N$ is as in  (9.3) with $\la_{ij}$ being the only non-zero entry. Moreover, $a\tau_w a\inv=\tau_{a_iw}$ for $\tau_w\in\W$ with $w\in W_i$.

If $h_i$ the linear functional on $\a$ projecting onto the $i$-th component, this shows that $\a$ acts on $\n+\w=\n_\G$ (the Lie algebra of $N_\G$) with weights $h_i-h_j$ ($0\le j<i\le n$) and $h_i$ ($0\le i\le n$). None of these weights being trivial, it follows that $\a$ is maximal abelian in $\p$ and the weights just obtained form a system of positive restricted roots of $\a$.

Using the isomorphism 
$$
\gathered
\a=\{(t_0,\dots,t_n)\}\longrightarrow \{(\tilde t_0,\dots,\tilde t_{n+1}):\sum\tilde t_i=0\}\\
(t_0,\dots,t_n)\longmapsto \Big(t_0-\frac1{n+2}\sum t_i, \dots, t_n-\frac1{n+2}\sum t_i,-\frac1{n+2}\sum t_i\Big)\ ,
\endgathered\tag10.5
$$
all the roots take the typical form $\tilde h_i-\tilde h_j$ ($0\le i,j\le n+1$, $i\ne j$) of the $A_{n+1}$-system. The fact that this system is irreducible implies that $\g$ is simple.

The last part of the statement is now obvious.
\endproof

We recall that a smooth map between two Riemannian manifolds is called {\it conformal} if its differential at any point is a scalar multiple of an orthogonal transformation.

\proclaim{Theorem 10.4} The action of collineations on projective $C$-lines is conformal. 
\endproclaim

\proof The elements of $U$ are obviously conformal. Composing with elements of $U$, we see that it suffices to prove the statement at points in $W$.

Suppose that $g$ is a translation in $\W$. Let $w_0+Cw$ be a tangent $C$-line at a point $w_0$. By (6.3), if $\ph$ is the angle of $Cw$ and $Cw_0$, $|X|_{w_0+}=c_{|w_0|,\ph}|X|$ for all tangent vectors $X\in Cw$. Similarly, $|g_*X|_{g(w_0)+}=|X|_{g(w_0)+}=c_{|g(w_0)|,\psi}|X|$, where $\psi$ is the angle of $Cg(w_0)$ and $Cw$. This implies the statement for $g\in\W$.

If $g\in GL(W,C)$, Corollary 9.3 implies conformality on lines through 0. Composing with translations, it follows that $g$ is conformal on all lines in $W$. By the decomposition (10.1), this finishes the proof.
\endproof

Since projective $C$-lines are spheres of the same dimension as $C$ (that we are assuming strictly greater than one), by Liouville's theorem conformal transformations are expressed by fractional linear transformations. For $C=\R$, Theorem 10.4 can be reformulated by replacing the word ``conformal'' with ``fractional linear''.
This also suggest that, when $V$ is trivial, the group $\Cal G$ that naturally replaces the collineation group is the conformal group when $C\ne\R$, and the group of fractional linear transformations  when $C=\R$. 

\vskip.2cm

Finally, the collineation group is the group of ``basic transformations'' of \cite{Te}, characterized (when $CPW$ is not a sphere) by the property of mapping Helgason spheres into Helgason spheres.

\vskip.5cm

\head 11. Compact rank-one spaces and symmetric cones \endhead

\vskip.3cm

We show that the groups $\G$ are exactly the groups such that $\G\times\R^+$ is the automorphism group of an irreducible symmetric cone \cite{FK}. In particular, $\G/U$ can be imbedded as a domain in real projective space in such a way that $\G$ acts by projective maps.

Taking into account Theorem 10.3, our statement will follow from the following  theorem, whose proof is based on representation theory.

\proclaim{Theorem 11.1} Let $\g$ be a real simple  non-compact Lie algebra whose restricted root system is of type $A_\ell$ and let $G$ be the adjoint group. Then $\tilde G=G\times \R^+$ is the automorphism group of an irreducible symmetric cone.
\endproclaim

The converse is well known, e.g. \cite{FK}, p. 108.

\proof
We write $\g=\k+\p$ for the Cartan decomposition and denote by $\theta$ the Cartan involution. We choose $\a\subset\p$ maximal abelian, and we complete its complexification $\a^\C$ to a $\theta$-invariant Cartan subalgebra $\h^\C$ of $\g^\C$. We identify $\h^\C$ with its own dual under the Killing form $(\,\cdot\,|\,\cdot\,)$ and write $\h$ for the real span of the roots. We also consider $\tilde\g=\g+\R$, $\tilde\a=\a+\R$, etc., and extend $(\,\cdot\,|\,\cdot\,)$ to an inner product on $\tilde\h$ so that $\R$ is orthogonal to $\h$. The theory of roots, restricted roots, weights, etc. still applies to the reductive algebra $\tilde\g$. If $\mu$ is a weight on $\tilde\h$, we denote by $\bar\mu$ its restriction to $\tilde\a$.

We choose an orthonormal basis $\{\eps_1,\dots,\eps_{\ell+1}\}$ of $\tilde\a\subset \tilde\h$ so that the restricted roots of $\tilde \g$ are $\bar\eps_i-\bar\eps_j$ ($i\ne j$). Since the Weyl group is transitive on the roots, the corresponding root spaces $\g^{\bar\eps_i-\bar\eps_j}$ have the same dimension, to be denoted by $d$. The $\R$-part of $\tilde\a$ is spanned by $\sum\eps_i$. We fix $\bar\eps_i-\bar\eps_j$ with $i< j$ as the positive roots.

Let $(\rho, V)$ be the irreducible representation of $\tilde \g$ with highest weight $2\eps_1$. By Theorem 4.12 in \cite{Hel3}, such a representation exists and it is a spherical representation, i.e. $V$ has a $K$-fixed vector. Indeed, $2\eps_1$ is in $\tilde\a$ and it satisfies the integrality condition
$$
\frac{(2\eps_1|\eps_i-\eps_j)}{(\eps_i-\eps_j|\eps_i-\eps_j)}\in\Z\qquad (\forall\,i\ne j)\ .
$$

We claim that $\dim_\C V=\dim_\R\,\tilde G/K$. To see this, we compute the dimensions of the $\a$-weight spaces $\bar V_{\bar\mu}=\{v\in V:\rho(H)v=\bar\mu(H)v\,,\,\forall\,H\in\a\}$. Each $\bar V_{\bar\mu}$ is the sum of the $\h$-weight spaces $V_{\mu'}=\{v\in V:\rho(H)v=\mu'(H)v\,,\,\forall\,H\in\h\}$, for those $\h$-weights $\mu'$ whose restriction to $\a$ is $\bar\mu$. 

By the general theory (cf. \cite{Hum}, p. 108, 114), the $\h$-weights $\mu$ of $\rho$ have the form $2\eps_1-\al$, where $\al$ is a positive $\h$-root and $(\al|2\eps_1)>0$. It follows that the only $\h$-weight of $\rho$ which restricts to $2\bar\eps_1$ is $2\eps_1$. Hence $\bar V_{2\bar\eps_1}=V_{2\eps_1}$ is one-dimensional.

By action of the Weyl group, $\bar V_{2\bar\eps_j}$ is one-dimensional for $1\le j\le \ell+1$. The other $\a$-weights are contained in the convex hull of the $2\bar\eps_j$ and equal to $2\bar\eps_j$ minus a sum of simple restricted roots. Therefore the other $\a$-weights can only be $\bar\eps_i+\bar\eps_j$, with $i<j$. 
By invariance under the Weyl group again, they must have the same multiplicity, so that we can restrict our attention to the weight $\bar\eps_1+\bar\eps_2$.

It is clear that $\bar V_{\bar\eps_1+\bar\eps_2}$ is the direct sum of all $V_{2\eps_1-\al}$ with $\bar\al=\bar\eps_1-\bar\eps_2$. We will prove that each such $V_{2\eps_1-\al}$ equals $\rho(\g_{-\al})V_{2\eps_1}$, with $\g_{-\al}$ the $\h$-root space relative to $-\al$, hence is one-dimensional. This will then show that $\dim\bar V_{\bar\eps_1+\bar\eps_2}=d$.

By \cite{Hum}, p.108, $V_{2\eps_1-\al}$ is spanned by the vectors
$$
v=\rho(X_{-\al_k})\cdots\rho(X_{-\al_1})v^+, \tag11.1
$$
where $v^+$ is a fixed non-zero element in $V_{2\eps_1}$, the $\al_j$ are positive $\h$-roots, $\sum_{i=1}^k\al_i=\al$ and $X_{-\al_i}$ spans the $\h$-root space $\g_{-\al_i}$. We may assume that $v\ne0$.

Since each $\bar\al_i$ is a positive restricted root or zero, there is one $i_0$ for which $\bar\al_{i_0}=\bar\eps_1-\bar\eps_2$, while $\bar\al_i=0$ for $i\ne i_0$. If $\bar\al_1=0$, then $2\eps_1-\al_1$ is not a weight for $\rho$, hence $\rho(X_{-\al_1})v^+=0$, a contradiction. So $i_0=1$. Now we prove that, if $k>1$, the number of factors in (11.1) can be reduced by one. Indeed, since $\rho(X_{-\al_2})v^+=0$, 
$$
\rho(X_{-\al_2})\rho(X_{-\al_1})v^+=\rho\big([X_{-\al_2},X_{-\al_1}]\big)v^+=\rho(X_{-\al'_1})v_+\ ,
$$
with some $X_{-\al'_1}\in\g_{-\al'_1}$, $\al'_1=\al_1+\al_2$.

Repeating this argument we find that $v$ is in $\rho(\g_{-\al})V_{2\eps_1}$.

By invariance under the Weyl group, $\dim_\C \bar V_{\bar\eps_i+\bar\eps_j}=d$ for all $i<j$, and adding up the dimensions, we find that $\dim_\C V=\ell+1+\frac{\ell(\ell+1)}2d$.
In $\tilde\g=\k+\tilde\p$ we have $\tilde\p=\tilde\a+\sum_{i<j}(\id-\theta)\g_{\bar\eps_i-\bar\eps_j}$. Counting the dimensions, we find that $\dim_\C\,V=\dim_\R\,\tilde G/K$.

We imbed $\tilde G/K$ into $V$ by choosing a $K$-invariant vector $e$ and defining the map $gK\mapsto \rho(g)e$. Then $V^0=\rho(\tilde\p)e$ is the tangent space of $\Omega=\rho(\tilde G)e$ at $e$ (under the usual identification). $V^0$ is a real form of $V$ and it is $\rho(K)$-invariant. Introducing an inner product invariant under the compact form of $G^\C$, $\rho(\tilde \p)$ consists of Hermitian linear transformations. Then $\rho(\tilde \a), \rho(\exp\tilde\a)$ are simultaneously diagonalizable and real. It follows that $\Omega=\rho(K)\rho(\tilde A)e$ is in $V^0$. Since its tangent space at $e$ is all of $V^0$, it is open.
Since the $\R^+$-part of $\tilde G$ acts by positive scalar transformations, $\Omega$ is a cone in $V^0$.

To see that $\Omega$ is symmetric in the sense of \cite{FK}, i.e. self-dual, we observe that our choice of inner product guarantees that $\rho(\tilde G)$ is closed under taking adjoints. This is enough to prove that $\Omega$ is self-dual, cf.  \cite{FK}, p. 20, Exercise 8.
\endproof

One can also describe $\Omega$ in more detail. By the Cartan decomposition, $V^0=\rho(K)V^{\rad}$, where $V^{\rad}$ is the subspace $\rho(\tilde\a)e$. Properly normalizing the inner product, the vectors $e_i=\rho(\eps_i)e$ form an orthonormal basis of $V^{\rad}$. Then $\Omega^{\rad}=\Omega\cap V^{\rad}$ is just the positive quadrant in $V^{\rad}$, and $\Omega=\rho(K)\Omega^{\rad}$.

We also note that the symmetric space $G/K$ can be realized as $\rho(K)\rho(\exp_G\a)e$, where $\rho (\exp_G\a)e$ is the hyperboloid $\{\sum_1^{\ell+1}t_ie_i: t_1t_2\cdots t_{\ell+1}=1\}$ in $V^{\rad}$.

Being a dual cone, $\Omega$ is convex. Its extremal generators can only be the $\rho(K)$-images of the edges in $\Omega^{\rad}$ and $K$ acts transitively on these. So $\E_1$, the intersection of the set of extremal generators with the unit sphere, is just $K\cdot e_1$. The stabilizer of $e_1$ in $K$ is the same as the centralizer of $\eps_1$ in $K$, and this is the $M'$-part of the Langlands decomposition of the parabolic subgroup $P=M'A'N'$ determined by the simple roots $\eps_i-\eps_{i+1}$ with $2\le i\le \ell$.

We identify $e_1$ and $\E_1$ with their images $\tilde e_1,\tilde\E_1$ in the projective space $PV^0$ and write $\tilde\rho$ for the action of $G$ on $PV^0$ induced by $\rho$. We claim that $\tilde\rho(G)\tilde e_1=\tilde \E_1$ and the stabilizer of $\tilde e_1$ is $P$. 
For this, we consider $\rho(\exp\, t\eps_1)e=e^te_1+e_2+\cdots+e_{\ell+1}$, which is fixed under the conjugate $K^{\exp\,t\eps_1}$ of $K$. Passing to $PV^0$ and letting $t\to\infty$, $\tilde e_1$ is fixed under the limit of $K^{\exp\,t\eps_1}$, which is $M'N'$, and also under $A'$, which acts on $e_1$ by scalars. This shows that $\tilde \E_1$ is one of the Satake-F\"urstenberg boundaries of $G/K$.

It is a result of U. Hirzebruch (cf. \cite{FK}, p.78, Exercise 5) that, for any irreducible symmetric cone, $ \E_1$ is a compact symmetric space of rank one, and (by classification) every such space arises in this way. We can now reprove this result, and a little more, if we show that, starting with any compact rank-one symmetric space, realizing it as in Section 5, then using Theorem 11.1 to construct the corresponding cone and its $\E_1$, we get back to the initial space.

For this, it is enough to check that when we apply Theorem 11.1 with $\G$ and $U$ in place of $G$ and $K$, the group $M'$ will be the same as the group $K$ of Sections 4-10. Now, in the Lie algebra of it, we can take $\tilde h_{n+1}$ as the element corresponding to $\eps_1$. After the coordinate change (10.5), this is a scalar multiple of $h_0+\cdots+h_n$, which generates the group of scalar transformations of $W$. The stabilizer of this in $U$ is indeed $K$, by Corollary 6.5.

\vskip.5cm

\head Appendix. The non-compact symmetric spaces \endhead

\vskip.3cm

What follows is redoing \cite{CDKR2} from a different starting point, to make it more compatible with the approach to compact spaces that we have taken.
\vskip.3cm

Let $B_W$ be the open unit ball in $W$. The tangent space $T_wB_W$ at a point in $w\in B_W$ is naturally identified with $W$ itself.  We introduce a Riemannian metric on $B_W$ by assigning, on the tangent space $T_wB_W$ at $w\in B_W$, the scalar product $\lan \ ,\ \ran_{w-}$ such that
$$
\lan X,Y\ran_{w-}=\left\{\matrix \displaystyle\frac{\lan X,Y\ran}{(1-|w|^2)^2}\hfill&\text{ if }X,Y\in
Cw\ ,\hfill\\ 
\displaystyle\frac{\lan X,Y\ran}{1-|w|^2}\hfill&\text{ if }X,Y\in
(Cw)^\perp\ ,\hfill\\
\displaystyle0\hfill& \text{ if } X\in Cw, Y\in (Cw)^\perp\ ,\hfill
\endmatrix\right.\tag A.1
$$
for $w\ne0$ and, passing to the limit for $w\to0$, 
$$
\lan X,Y\ran_{0-}=\lan X,Y\ran\qquad \text{ for every }X,Y\ .\tag A.2
$$

In \cite{CDKR2} these formulas are obtained towards the end, in Section 7; here they are definitions. Note that in \cite{CDKR2} this metric is multiplied by a factor 4. This multiplies by a constant the arc lengths and curvatures, but does not change anything essential such as isometries or geodesic submanifolds. 

It is immediate that the same $K$ of Section 4 consists of isometries fixing the origin.
Adapting the proofs in Section 6, one easily obtains the following properties.

\proclaim{Lemma A.1} The geodesic $\gamma$ going through the origin with tangent vector $w\in S_W$ is the diameter in the direction of $w$, parametrized as $\gamma(t)=(\tanh t) w$.
 The group of isometries of $B_W$ fixing the origin is $K$.
\endproclaim

The maps $a_t$ defined by
$$ 
a_t(\zeta,v)=
\big((\sinh t\,\zeta+\cosh t\one)\inv(\cosh t\,\zeta+\sinh t\one),
(\sinh t\,\zeta+\cosh t\one)\inv v\big)\ .\tag A.3
$$
form a one-parameter group $A$ of isometries of $B_W$. This can be proved by a computation very similar to the proof of Lemma 6.2.

Since the orbit of $(0,0)$ under $A$ consists of the points $a_t(0,0)=(\tanh t,0)$, it follows that the group $G$ of isometries of $B_W$ is transitive on $B_W$. We also see that $G=KAK$. Since the isometry $w\mapsto -w$ has a single fixed point, it also follows that $B_W$ is a symmetric space. Transitivity of $K$ on $S_W$ implies that the rank of $B_W$ is one.

Simple modifications to the proof of Theorem 7.2 give that the sectional curvature of a plane element $E$ spanned by $X,Y\in W$ in the tangent space to $B_W$ at the origin is the negative of $\sigma(E)$ in (7.1). 

The unit balls $B_1$ in $W_1=C_1\oplus V_1$ and $B_2$ in $W_2=C_2\oplus V_2$ are isometric if and only if $(C_1,V_1)\sim(C_2,V_2)$ as $C$-modules. However, the degenerate cases $(\R,V)$ and $(\R\oplus V,0)$ give different models of the same space. The map
$$
(t+v,0)\longmapsto \Big(\frac{2t}{1+t^2+|v|^2},\frac{2v}{1+t^2+|v|^2}\Big)
$$
from the second to the first is an isometry up to a factor 2. For $V\ne0$, these are respectively the Klein model and the Poincar\'e model of real hyperbolic space. For $V=0$, they trivially become models of the real line.
\cite{CDKR2} contains a classification-independent proof of the fact that the metrics (A.1) give all the rank-one symmetric spaces of the non-compact type (together with $\R$, the irreducible Euclidean non-compact symmetric space).

The maps $a_t$ in (A.3) extend to homeomorphisms of $CPW$. The extension to $W\setminus\big\{(\zeta,v):\zeta\ne -\coth t\big\}$ is obvious, and the extension to the rest of $CPW$ is, in analogy with (6.10),
$$
\align
a_t(-\coth t,v)&=[\one,-(\sinh t)v]\ ,\\
a_t\big([\one,v]\big)&=\big(\coth t,(\sinh t)\inv v\big)\\ a_t\big([0,v]\big)&=[0,v]\ .
\endalign
$$

Therefore, all isometries of $B_W$ extend to homeomorphisms of $CPW$. In analogy with Corollary 7.4, it follows that these extensions are collineations.

\vskip.3cm

To tie up with \cite{CDKR2}, we use the Cayley transform $c$ (cf. \cite{CDKR2}, p. 208), which in our notation is 
$$
(\zeta',v')=c(\zeta,v)=\big((\one-\zeta)\inv(\one+\zeta),2(\one -\zeta)\inv v\big)\ .
$$

Clearly, 
$$
\RE\zeta'-\frac14|v'|^2=\frac{1-|\zeta|^2-|v|^2}{|\one-\zeta|^2}\ ,
$$
which shows that $c$ maps $B_W$ onto 
$$
D=\big\{(\zeta',v'): \RE\zeta'-\frac14|v'|^2>0\big\}\ .
$$

Define $B:V\times V\rightarrow C$ by
$$
\lan B(v,v'),\zeta\ran=\lan \zeta v,v'\ran\ .
$$

Among the isometries of $D$, with the metric induced on it by $c$, one has 
 the ``translation'' group $\tilde N=\{\tilde n_{(z,u)}:z\in
C'\,,\,u\in V\}$, with
$$
\tilde n_{(z,u)}(\zeta,v)=\Big(\zeta+z+ \half B(v,u)+\frac14|u|^2,v+u\Big)\ ,
$$
acting simply transitively on the level sets of the height function
$h(\zeta,v)=\RE \zeta-\frac14|v|^2$.

The Lie algebra $\n$ of $\tilde N$ is an $H$-type Lie algebra satisfying the $J^2$-condition. The discussion in Section 3 shows that all such algebras show up in this way.

If $\tilde a_t=ca_tc\inv$, then 
$$
\tilde a_t(\zeta,v)=(e^{2t}\zeta, e^tv)\ ,
$$
the group $\tilde A=cAc\inv$ normalizes $\tilde N$, and $\tilde N\tilde A$ acts simply transitively on $D$. This is the starting point in \cite{CDKR2}. 

We note that $c$ is almost the same as $b_{\frac\pi4}$ defined in Lemma 6.2, except that the factor 2 in $c$ is replaced by $\sqrt2$ in $b_{\frac\pi4}$. One can in fact replace $c$ by $b_{\frac\pi4}$ and $D$ by 
$D'=\big\{(\zeta',v'): \RE\zeta'-\frac12|v'|^2>0\big\}$, a choice that may be preferable in some ways.

\enddocument